\documentclass[onefignum, onetabnum]{siamonline220329}

\usepackage[T1]{fontenc}
\usepackage[utf8]{inputenc}
\usepackage{amsmath}
\usepackage{amssymb}
\usepackage{verbatim}
\usepackage{color}
\usepackage{array}
\usepackage{graphicx}
\usepackage{epstopdf}
\usepackage{upquote}
\usepackage{epsfig}
\usepackage{mathtools}
\usepackage{courier} 
\usepackage{xcolor} 
\usepackage{bm} 

\usepackage[]{hyperref}
            
\newcommand{\be}{\begin{equation}}
\newcommand{\ee}{\end{equation}}

\newcommand{\red}[1]{{\color{black}{#1}}}

\usepackage[normalem]{ulem}

\usepackage{lipsum}
\usepackage{amsfonts}
\usepackage{graphicx}
\usepackage{epstopdf}
\usepackage{algorithmic}
\ifpdf
  \DeclareGraphicsExtensions{.eps,.pdf,.png,.jpg}
\else
  \DeclareGraphicsExtensions{.eps}
\fi

\newsiamremark{remark}{Remark}
\newsiamremark{hypothesis}{Hypothesis}
\crefname{hypothesis}{Hypothesis}{Hypotheses}
\newsiamthm{claim}{Claim}

\headers{Approximate solutions to a nonlinear functional differential equation}{N. Hale, E. Thomann, JAC. Weideman}

\title{Approximate solutions to a nonlinear functional differential equation }

\author{Nicholas Hale\footnote{Corresponding\;author.}\,\;\thanks{Dept of Mathematical Sciences, Stellenbosch University, Stellenbosch, 7600,  South Africa (nickhale@sun.ac.za)}
\and Enrique Thomann\thanks{Dept of Mathematics, Oregon State University, Corvallis, OR 97331,  USA (enrique.thomann@oregonstate.edu)}
\and 
JAC Weideman\thanks{Dept of Mathematical Sciences, Stellenbosch University, Stellenbosch, 7600,  South Africa (weideman@sun.ac.za)}
}

\usepackage{amsopn}

\ifpdf
\hypersetup{
  pdftitle={Approximate solutions to a nonlinear functional differential equation},
  pdfauthor={N. Hale, E. Thomann, JAC. Weideman}
}
\fi

\begin{document} 
  


\maketitle  

\begin{abstract}
A functional differential equation related
to the logistic equation  
is studied by a combination of numerical and perturbation
methods.  Parameter regions are identified where
the solution to the nonlinear problem is approximated
well by known series solutions of the linear
version of the equation.   The solution space
for a {particular} class of functions is then mapped out 
using a continuation approach.
\end{abstract}




\section{Introduction}\label{sec:intro}%
 
The functional differential equation
\be
u^{\prime}(t) +a \, u(t) = b \, u(\alpha t), \quad t > 0, 
\label{eq:linear}
\ee
with $a$, $b$, and $\alpha > 0$ constants, was the topic of two major papers,~\cite{Fox} and~\cite{Kato}, both dating back about half a century.  An application 
mentioned in the former of these references is the dynamics of an overhead current collection 
system for an electric locomotive.  Even
earlier papers derived the equation as a probabilistic model 
of absorption processes~\cite{Ambartsumian,Gaver}.  In later years the equation has re-appeared
in mathematical biology, such as in the modelling of cell growth~\cite{Hall,VanBrunt}.  

In this paper we present results for a nonlinear version of this equation, 
\be
u^{\prime}(t) +  u(t) =  u^2(\alpha t), \quad t > 0, 
\label{eq:nonlinear}
\ee
considered by~\cite{Athreya} 
and more recently by~\cite{Dascaliuc2019}.  The latter
paper, in which the  term {\em $\alpha$-Riccati equation} is coined, discusses
an application to cascading processes in Fourier space representations 
of some nonlinear PDEs (including the Navier--Stokes equations).  
The two references cited both approach the nonlinear problem from
a probabilistic viewpoint, in particular a stochastic branching process. 
Here, the approach is numerical computation
supported by a perturbation analysis.
 
For $\alpha = 1$, both~(\ref{eq:linear}) and~(\ref{eq:nonlinear}) reduce to well-known equations that will not be considered here.   
  For $\alpha < 1$, the equations are delay differential equations  with proportional or `pantograph-type' delay. This regime was investigated  
  in~\cite{Fox} and~\cite{Kato} in the linear case, and~\cite{Dascaliuc2019} in the nonlinear
  case.  In particular,  existence and uniqueness of solutions were established
for both~(\ref{eq:linear}) and~(\ref{eq:nonlinear}).  

{The case of $\alpha>1$ in~(\ref{eq:linear}) and~(\ref{eq:nonlinear}) involves a `time-advanced' argument. However, in applications this {typically} represents a non-temporal variable, thus avoiding a violation of causality in the usual time-evolution sense.   For example, the independent variable {may} represent a scale in the biological settings treated in~\cite{Hall}, the magnitude of the wave vector in the self-similar Navier--Stokes equations as discussed in~\cite{Dascaliuc2019}, or a transparency scale as in~\cite{Ambartsumian}.  In the pantograph-type systems introduced in~\cite{Fox}, and carefully analysed in~\cite{Kato}, the parameter $\alpha$ is also allowed to take values greater than unity for modelling purposes.  Further examples in which $\alpha>1$ arises are models of nerve axons from conduction theory~\cite{chi1986} and preview control from control theory, where it {models} anticipation~\cite{birla2015}. 
As the simplest {nontrivial} nonlinear equation with this feature, we investigate solutions to~(\ref{eq:nonlinear}) in order to provide insight into more complex time-advanced functional differential equations arising in applications.

While~\cite{Fox} and~\cite{Kato} give a thorough analysis of properties of solutions of the linear equation~(\ref{eq:linear}), current knowledge for the nonlinear counterpart~(\ref{eq:nonlinear}) is incomplete.  A {visual} status report {is given in}~\cite[Figure~2]{Dascaliuc2019}.  This diagram catalogues the existence and uniqueness properties of solutions in several regions of the $(\alpha, u_0)$ plane, where $u(0)=u_0$ represents an initial condition. In several regions either existence or uniqueness remains unknown.  The focus here will be on solutions to~(\ref{eq:nonlinear})  for the case $\alpha>1,$ $u_0=1$, which forms a small but important section of that diagram. }

We shall  make no attempt to prove any 
theorems regarding existence or uniqueness, but rather focus on quantitative
aspects such as the
approximation of solutions, or, since multiple solutions exist, to 
determine how many solutions in a certain class are possible for each value of $\alpha$. 
By a combination of analysis and
computation we identify parameter regions where   
the solutions to the nonlinear problem (\ref{eq:nonlinear}) can be approximated well by 
certain explicitly known series solutions of the linear problem~(\ref{eq:linear}).

For numerical solutions
we use a spectral collocation method for functional differential equations. (Details  can be found in~\cite{Hale2023} and a short description in section~\ref{sec:nummeth}.)  This method provides more efficient 
and more accurate results than other methods reported in the 
literature.   Two such methods for the 
linear problem~(\ref{eq:linear}) are described in~\cite{Fox}:
A finite difference method, applicable to the case $\alpha<1$,
and a Lanczos tau-method for both $\alpha<1$ and $\alpha>1$.
Numerical results, of low accuracy in today's terms but
impressive for 1970s technology, are given. Because the problem~(\ref{eq:nonlinear}) is nonlinear, any
numerical procedure would require iteration and this demands
good starting guesses.   Such approximations, 
which may
be of interest in their own right,
are obtained here by an elementary perturbation analysis.
 
The only numerical approach to the
nonlinear problem~(\ref{eq:nonlinear}) that we are aware of is
the probabilistic simulation of~\cite{DascaliucErratum}.
Like most methods of Monte Carlo-type, however, it is not 
efficient for doing the extensive experimentation  that will be 
reported here. 

The outline of this paper is therefore as follows. In section~\ref{sec:linear} we introduce some aspects of the linear problem~(\ref{eq:linear}) that are pertinent in the analysis of the nonlinear problem~(\ref{eq:nonlinear}). In particular, we consider $a=1$ and $b = 2$, corresponding to the linearisation of~(\ref{eq:nonlinear}) about the constant solution $u = 1$, and demonstrate the existence of nontrivial solutions to the homogeneous linear problem for certain characteristic values of $\alpha$, designated here as  $\alpha_n$, $n = 1, 2, \ldots$. In section~\ref{sec:nonlinear} we show that these characteristic solutions, appropriately scaled, provide the leading order term in a perturbation of the constant solution to~(\ref{eq:nonlinear}) in the neighbourhood of these characteristic values. In section~\ref{sec:nummeth} we use these perturbation approximations in combination with the spectral collocation method and pseudo-arclength continuation to determine 
numerically a family of solutions to~(\ref{eq:nonlinear}) for a range of $\alpha$ values, which we discuss in section~\ref{sec:family}.

\section{Summary of results for the linear equation}\label{sec:linear}%

The papers~\cite{Fox} and~\cite{Kato} contain a wealth of information on the linear
equation~(\ref{eq:linear}).   This section highlights only the major
results that will be used in the nonlinear analysis of the next section. 
In particular, in that section the values $a=1$ and $b = 2$ will be relevant, i.e.,
\be\label{eq:linear2}
u^{\prime}(t) + u(t) = 2 u(\alpha t).
\ee
{A series solution can be obtained from~\cite[Thm 9]{Kato}, namely:}
\be
u(t) = CE(t; \alpha), \quad \text{where} \quad E(t; \alpha) =  e^{-t} \Big( 1 + \sum_{k=1}^{\infty}   
\frac{2^k e^{(1-\alpha^k)t}}{(1-\alpha)(1-\alpha^2)\cdots(1-\alpha^k)} \Big),
\label{eq:series}
\ee
and $C$ is a constant.\footnote{Neither~\cite{Kato}
nor~\cite{Fox} explicitly derive~(\ref{eq:series}) 
but~\cite{Hall} and~\cite{VanBrunt} contain derivations based on Laplace transforms and Mellin transforms, respectively.}  The series
converges absolutely and uniformly for all $t \geq 0$ in the case 
$\alpha>1$.

The initial condition $u(0)=u_0$ can be matched  to~(\ref{eq:series}) by setting
$t=0$, but of particular interest is the case $u_0=0$.   This can be satisfied
by setting $C=0$, but 
nontrivial solutions are possible for 
special values of $\alpha$, as the following lemma shows.  These 
values of $\alpha$ can therefore be
viewed as eigenvalues or characteristic values of the linear operator. 
In other words, for these values of $\alpha$ the function~(\ref{eq:series})
lies in the nullspace of~(\ref{eq:linear2}).  

\smallskip

{\sc Lemma 1:}  The equation 
\be
1 + \sum_{k=1}^{\infty}   
\frac{2^k}{(1-\alpha)(1-\alpha^2)\cdots(1-\alpha^k)}  = 0 
\ee
is satisfied if and only if $\alpha = \alpha_n = 2^{1/n}$ for some
positive integer $n$.

\smallskip

{\em Proof:}  An identity  well-known in the field
of $q$-series~\cite[p.~11]{Gasper} states that
\be
1 + \sum_{k=1}^{\infty} \frac{(-1)^k q^{k(k+1)/2} c^k}{(1-q)(1-q^2)\cdots (1-q^k)} = (1-c q)(1- cq^2) \cdots (1-c q^n) \cdots , \quad |q|<1.
\label{eq:euler}
\ee
The proof is completed by substituting $c = 2$, $q = 1/\alpha$. \qed

\smallskip

A related observation was made in~\cite{Kato},  namely that the initial condition $u(0) = 0$ ``can occur only for exceptional pairs $a$ and $b$''  (with these constants as defined in~(\ref{eq:linear})).  The
authors did not explore further, but in~\cite[Eq.~(2.13)]{VanBrunt} it was shown
that $u(0) = 0$   if and only if $b = a \alpha^n$ for some integer $n$.
This follows from~(\ref{eq:euler}) 
by substituting $c = b/a$ and $q=1/\alpha$.   Note that  both~\cite{Kato} and~\cite{VanBrunt} addressed the situation with
$\alpha$ fixed; our interest is in the case when $a$ and $b$ are fixed.

\begin{figure}[htb]
\begin{center}
{\hspace*{-10pt}\includegraphics[width=.95\textwidth] {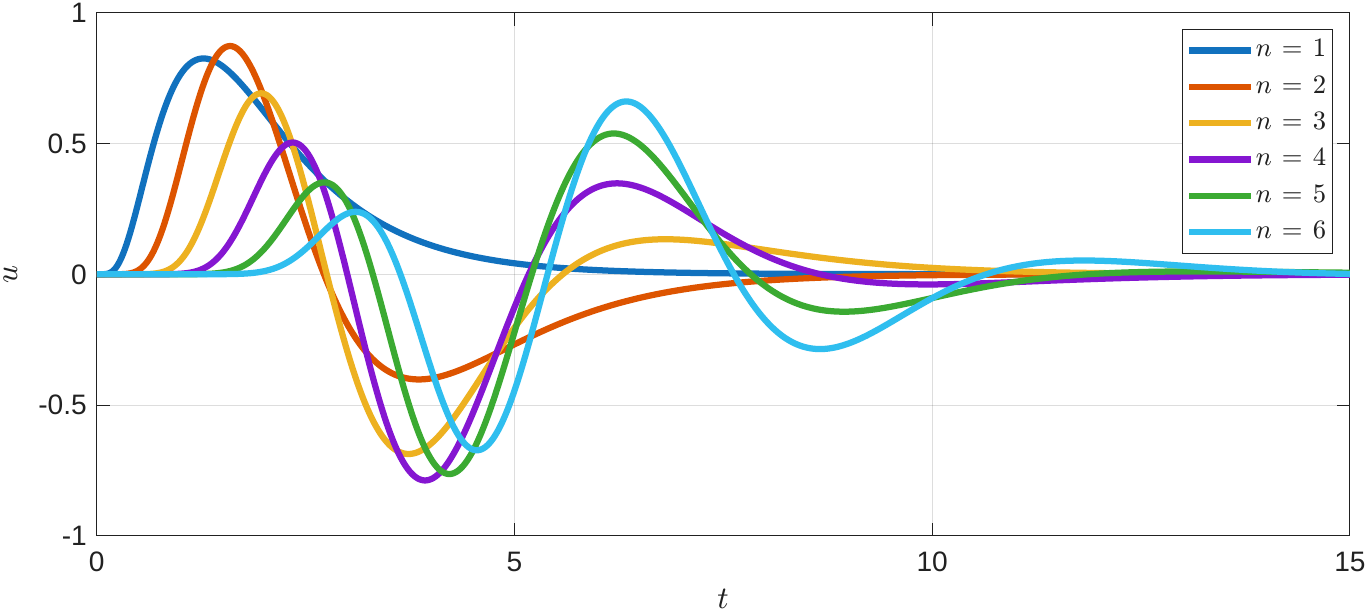}}
\caption{$L^2$ normalised characteristic solutions to~(\ref{eq:linear2}),  
as defined by the series (\ref{eq:series}),
for $\alpha_n = 2^{1/n}$, $n = 1, 2, \ldots, 6$. An additional scaling of $(-1)^{n+1}$ is also included so that all solutions are initially positive. 
}
    \label{fig:char}
\end{center}   
\end{figure}

To distinguish these cases we shall call the functions~(\ref{eq:series})
with $\alpha = \alpha_n = 2^{1/n}$ ``characteristic functions'' rather
than ``eigenfunctions'' as in~\cite{VanBrunt}.  The first few characteristic
functions of~(\ref{eq:linear2}) are shown in Figure~\ref{fig:char}, normalised such that $\| u \|_2 = 1$ (the norm is the standard, unweighted $L^2$ norm on 
$[0,\infty)$.) 

Noteworthy features of the characteristic solutions shown in Figure~\ref{fig:char} include
the observation that for $n=1$ the function has no zero in $(0,\infty)$, for $n = 2$ there is one zero, 
for $n = 3$ two zeros, and so on.   This is reminiscent of  Sturm--Liouville theory, but here the roots of two successive
characteristic functions do not strictly
interlace.    (See also~\cite{VanBrunt2} for 
comments about orthogonality, but keep in mind that their 
definition of an eigenfunction is different, as mentioned below Lemma~1.)
Another striking feature is the very flat profile near
$t=0$.  By repeated differentiation of~(\ref{eq:linear2}) one can show that if $u(0)=0$, then all derivatives of $u$ vanish at $t=0$, suggestive of an essential singularity
at the origin. On the other hand, Fox et al.~\cite[Thm.~2]{Fox} show two types of asymptotic behaviour as  $t\rightarrow\infty$, namely $u=O(e^{-t})$ {or  $u=O(t^{-\nu})$, $\nu = \log{2}/\log{\alpha}$.  The solutions shown in Figure~\ref{fig:char} are of the former type.}

\section{Perturbation analysis of the nonlinear equation}\label{sec:nonlinear}%
As mentioned in the introduction, the focus of this paper is on the nonlinear problem~(\ref{eq:nonlinear}) with initial condition $u(0) = 1$.   
In~\cite[Thm.~4.1]{Dascaliuc2019} two families of global solutions have been
identified according to the asymptotic
limits (A) $u \to 1$ as $t \to \infty$ and (B) $u \to 0$ as $t \to \infty$.  Here we restrict ourselves to the cases where these two limits are approached exponentially, i.e.,
according to (A) $u = 1+ {\rm O}(e^{-t})$  and (B) $u = {\rm O}(e^{-t})$.
{As in the linear case,} algebraic approaches to these limits may also be possible~\cite{Dascaliuc2024},  
but this will require a  different procedure from the one described here. 

An obvious solution in class (A) is the constant function $u = 1$.   Depending
on the value of $\alpha$, however, there may be one or more non-constant
solutions.  When computing these solutions  two complications arise, namely,
nonlinearity, and non-uniqueness.   A successful computation will therefore
require nonlinear iteration starting from  a  good initial guess.   
The numerical iteration procedure is described in section~\ref{sec:nummeth}.
Here we present an elementary perturbation
analysis for generating good initial guesses,  but which may be of interest in its own right.

Because our focus is on solutions in class (A), we {consider} $v(t)=u(t)-1$ {with $v = O(e^{-t})$}.
Substituting into~(\ref{eq:nonlinear}) gives
\be
v^{\prime}(t) + v(t) = 2 v(\alpha t) + v^2(\alpha t), \quad v(0) = 0.
\label{eq:v}
\ee
With the intuition (and support from numerical experimentation) that the characteristic values from the previous section act as bifurcation points in the nonlinear 
equations~(\ref{eq:nonlinear})  and (\ref{eq:v}), we seek solutions in the neighbourhood of these points by considering $\alpha = \alpha_n+\epsilon$, $|\epsilon| \ll 1$. Further, assuming small perturbations about the trivial solution $v(t) = 0$, we try the {formal} regular 
expansion $v(t) = \epsilon v_0(t) + \epsilon^2 v_1(t) + \cdots$, {where  the $v_k(t)$ are assumed to be bounded on $t >  0$}.  This gives,
to leading order,  
\be
v_0^{\prime}(t) + v_0(t) = 2 v_0(\alpha_n t),  \quad v_0(0) = 0,
\ee
which is the linear equation discussed in section~\ref{sec:linear}.\footnote{In an effort to avoid overencumbering notation, we omit burdening $v_0$ with an additional subscript, $n$. However, it should be understood here, and throughout, that each $v_0$ corresponds to a particular $\alpha = \alpha_n$. This value of $n$ should be clear from context.} 
By  Lemma~1, 
the explicit solution
to this equation is given by $v_0 = CE(t;\alpha_n)$, where $E(t;\alpha_n)$ is   defined in~(\ref{eq:series}). The constant $C$ remains to be determined, however, and
this can be done by matching solutions of the linear 
and nonlinear problems.   

To do this, we {follow~\cite{Hall} and} consider the moments of $v$ as follows.  
For arbitrary $\alpha$,
integrate~(\ref{eq:v})  with respect  to $t$ to obtain
\be
\int_0^\infty v^{\prime}(t) \, dt + \int_0^\infty v(t) \, dt
= 2 \int_0^\infty v(\alpha t) \, dt + \int_0^\infty v^2(\alpha t) \, dt.
\label{eq:mom}
\ee
Under the assumption that $v(0) = 0$ and $v \to 0$ as $t \to \infty$, the first
term vanishes.   Making a change of variable $t\mapsto t/\alpha$  on the right, multiplying by
$\alpha$, and rearranging gives
\be
(\alpha-2) \int_0^\infty v(t) \, dt =  \int_0^\infty v^2( t) \, dt.
\label{eq:moment1}
\ee
Consider now the particular perturbation solution corresponding to $n=1$, so that   $\alpha = 
2 + \epsilon$.  Inserting this and $v(t) = \epsilon v_0(t) + \epsilon^2 v_1(t) + \cdots$ 
into~(\ref{eq:moment1}) and then collecting terms gives, at leading order, 
\be
\int_0^\infty v_0(t) \, dt =  \int_0^\infty v_0^2( t) \, dt.
\ee
By {substituting} the series solution $v_0 = C_1E(t;2)$ from~(\ref{eq:series}),  one   finds that   the required value
of $C_1$ is 
\be
C_1 = \frac{\displaystyle\int_0^\infty E(t; 2) \, dt}{\displaystyle\int_0^\infty E^2(t; 2) \, dt } = 11.36911520\ldots,
\label{eq:C1}
\ee%
where the numerical value has been computed by 
the formulas (\ref{eq:summations}) in the Appendix. 
Similar formulas for $C_n$, $n >1$, can be derived by considering 
higher moments of $v$.  Their numerical values are given in Table~\ref{table:c} but details of the derivation are deferred to the Appendix. 

\begin{table}[htbp!]
\begin{center}
\begin{minipage}[]{.48\textwidth}
\begin{tabular}{c|rrr}
\multicolumn{1}{>{\centering\arraybackslash}m{10mm}|}{$n$} 
    & \multicolumn{1}{>{\centering\arraybackslash}m{23mm}}{$\hspace*{16pt}C_n $} 
    & \multicolumn{1}{>{\centering\arraybackslash}m{23mm}}{$\hspace*{13pt}\|C_nE_n\|_2$}
    \\ \hline
     1  &     11.36911520 &    1.811978234 \\
     2  &  $-$809.3665721   &   9.266935159 \\
     3  &   31551.15567   &   29.58298033 \\
     4  &  $-$1099159.137  &   85.46096526 \\
     5  &   34825078.48 &   224.4734290 \\
     6  &  $-$1045480822.  &   557.2788228 
\end{tabular}
\end{minipage}
\end{center}
\vspace*{.5em}
\caption{The first six scaling coefficients, $C_1, \ldots, C_6$, as given by~(\ref{eq:scalingcoeffs}), accurate to all digits shown. Although the $C_n$ grow rapidly with $n$, {the norms of} the characteristic functions, $E_n(t) = E(t,\alpha_n)$, are at the same time decreasing in magnitude. Their product, which  is the coefficient of $\epsilon$  in
the perturbation approximation~(\ref{eq:pert}), also increases in magnitude, but not as rapidly as $C_n$. Nevertheless, this means that as $\alpha$ approaches the value $1$ with $\epsilon$ fixed, the perturbation
approximations become less accurate. 
} \label{table:c}
\end{table}

In summary, the approximation to~(\ref{eq:v}) corresponding to $\alpha = 2^{1/n}+\epsilon$
is given to first order by 
\be
v(t) \approx  C_n E(t; 2^{1/n}) \epsilon,
\label{eq:pert}
\ee%
where $E(t; 2^{1/n})$ is the series solution~(\ref{eq:series}) and the $C_n$ values are given 
in Table~\ref{table:c}.

\begin{figure}[t!]  
\centering
     \includegraphics[height=94pt]{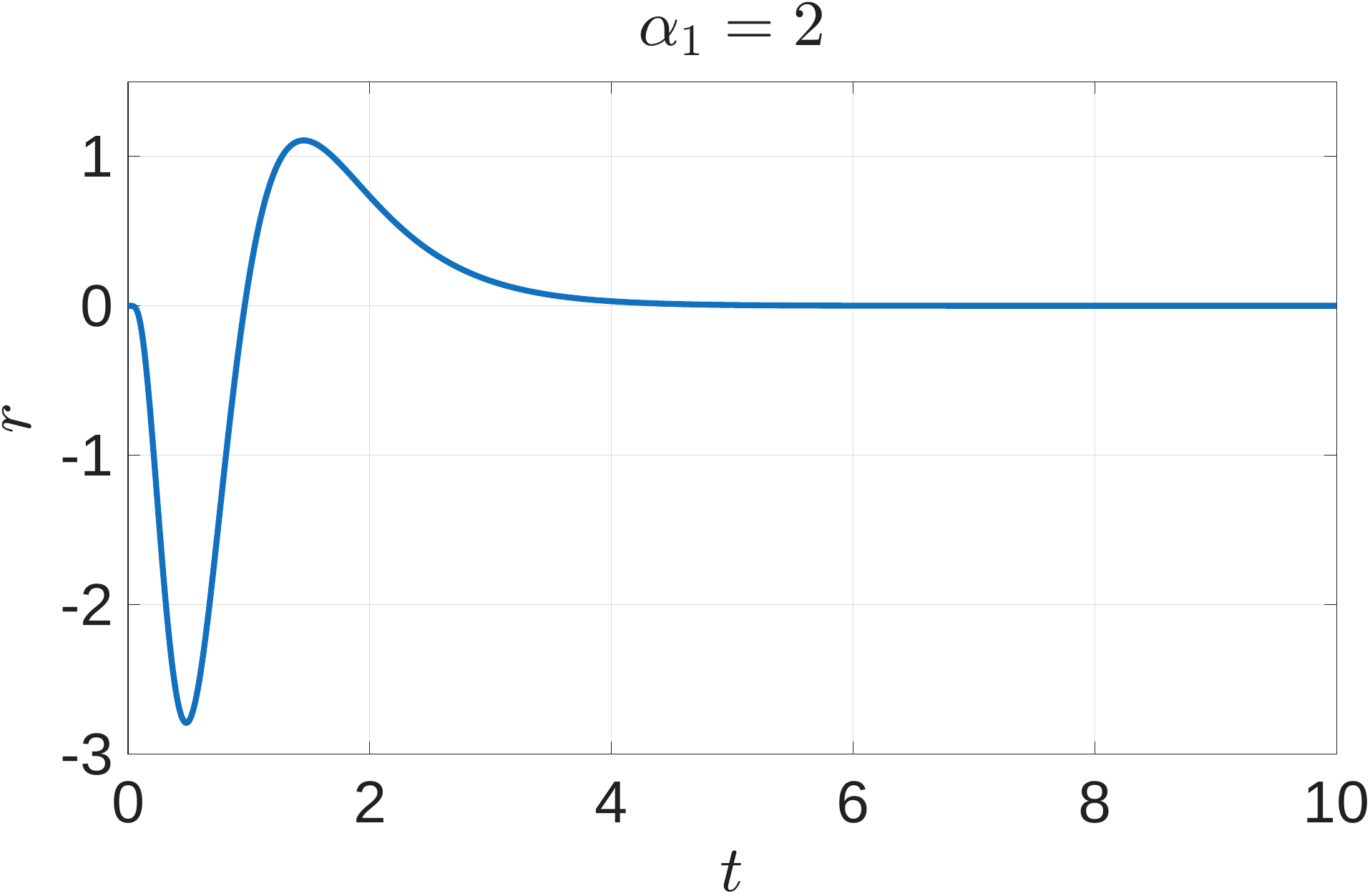}  %
    \includegraphics[height=94pt, trim={22pt 0 0 0},clip]{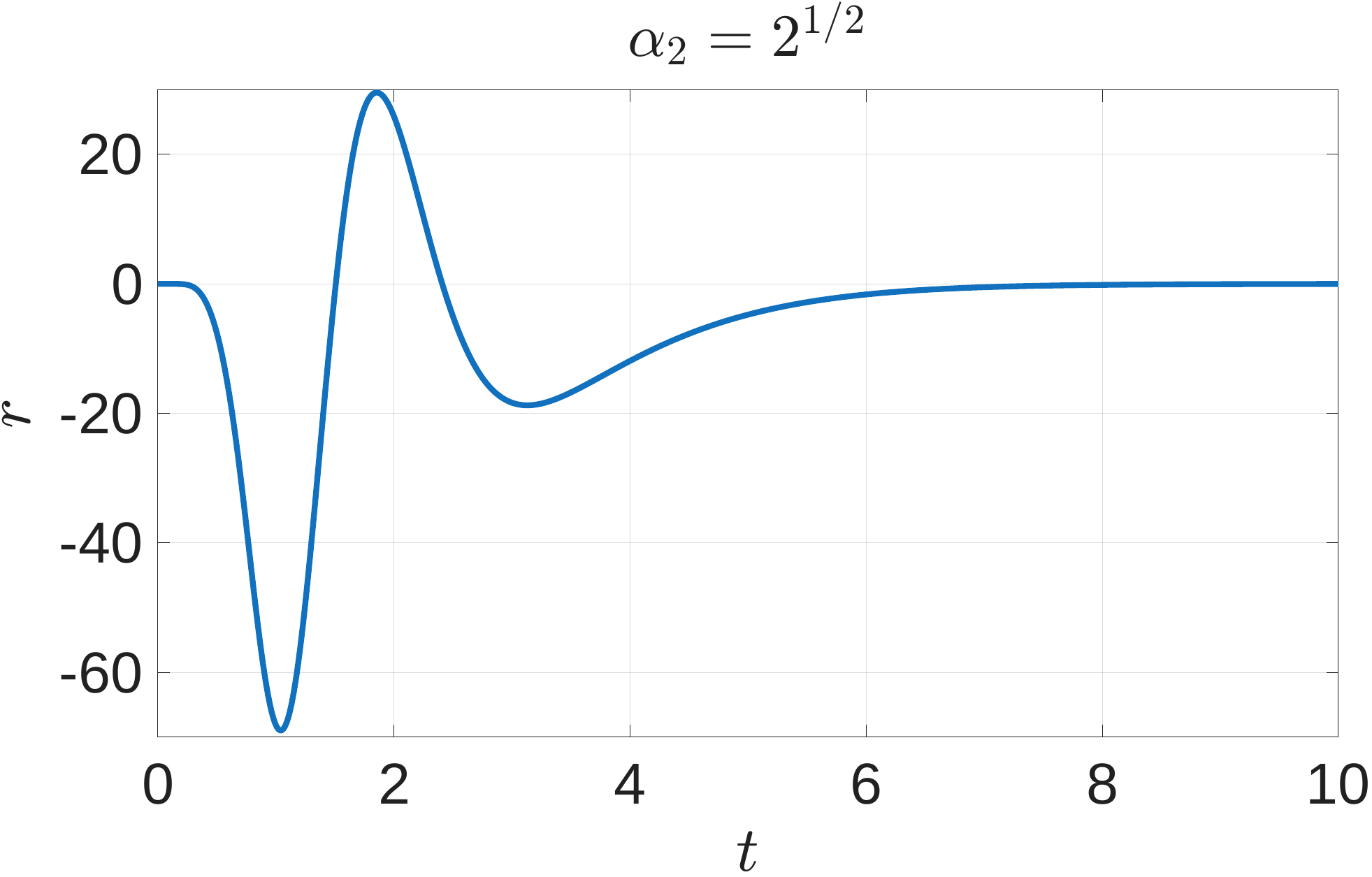}  %
     \includegraphics[height=94pt, trim={22pt 0 0 0},clip]{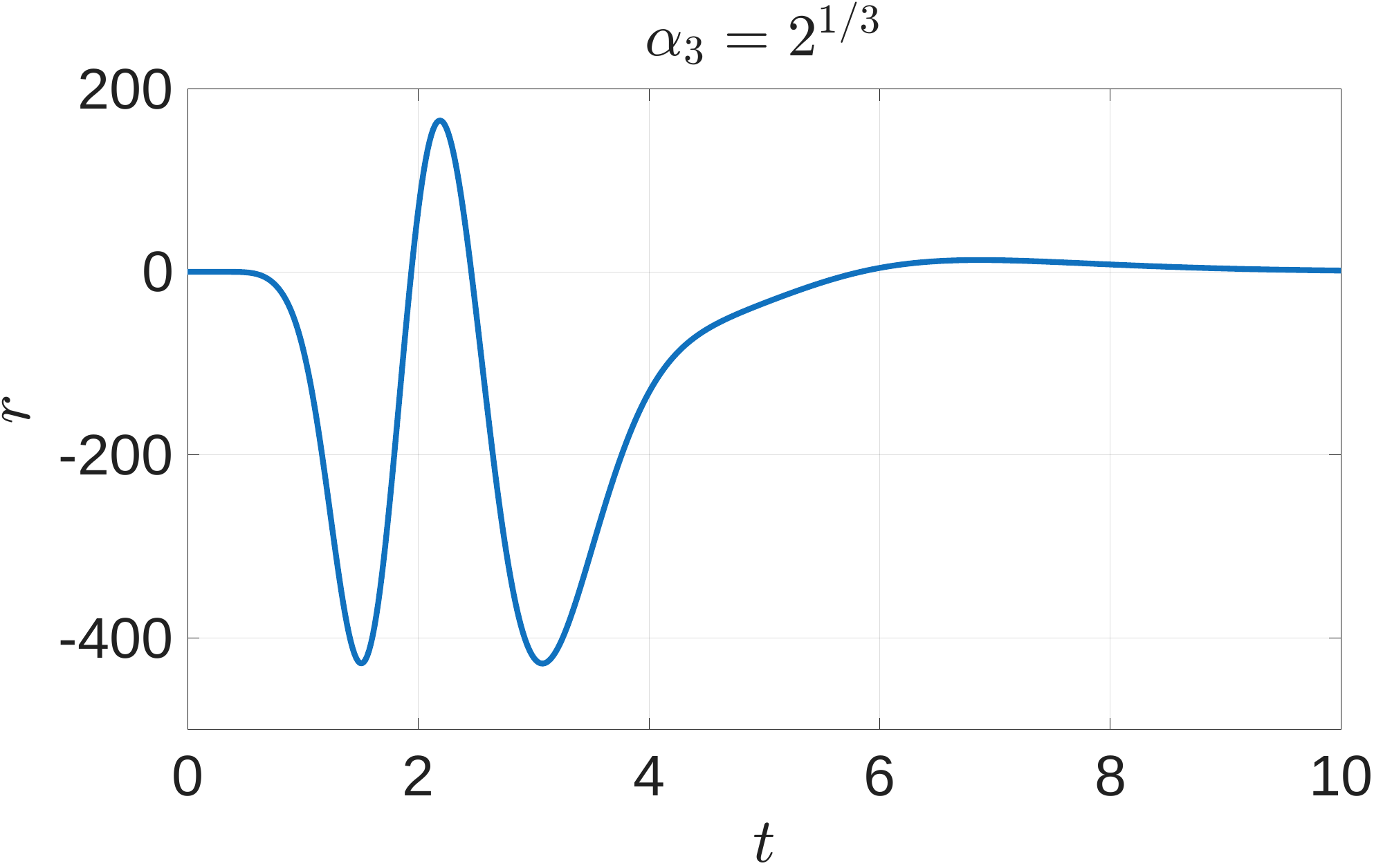}\\[-1em]
    \caption{The coefficient functions $r_n(t;\epsilon)$ as defined by~(\ref{eq:BE}), in the limit as $\epsilon \to 0$.   The observation that
    these functions are bounded on $[0,\infty)$ provides support for the
    validity of the formulas~(\ref{eq:pert}) as uniform approximations to the solutions of 
    (\ref{eq:v}). {These figures remain virtually unchanged for values of $|\epsilon|$ up to $10^{-2}$, indicating that the bounds hold not just in the limit, but in practice.}}
    \label{fig:BE}
\end{figure}

We have not pursued a formal proof of the validity of~(\ref{eq:pert}) in any
asymptotic sense.  (The analysis is rather complicated and perhaps of secondary
importance because the 
principal use of~(\ref{eq:pert}) is to provide initial guesses for the numerical method.) Instead, here is a computational check.  
We substitute~(\ref{eq:pert}) into~(\ref{eq:v}), 
and,  expecting the remainder to be ${\rm O}(\epsilon^2)$, we write
\be
v^{\prime}(t) + v(t) - 2 v((2^{1/n}+\epsilon) t) - v^2((2^{1/n}+\epsilon) t) = r_n(t; \epsilon) \, \epsilon^2,
\label{eq:BE}
\ee
for some coefficient function, $r_n$.
Every quantity on the left can be explicitly computed from the series~(\ref{eq:series})
and the tabulated values of $C_n$. Therefore the values of $r_n(t;\epsilon)$ 
can be computed for each $\epsilon$.  By Taylor analysis,
the limiting expressions 
as $\epsilon \to 0$ are
found to be 
\be
r_n(t;0) = -2C_n E^{\prime}(2^{1/n} t; 2^{1/n}) t - \big(C_n E(2^{1/n}t; 2^{1/n})\big)^2.
\label{eq:rn}
\ee
These functions   are shown for a few values of $n$ in Figure~\ref{fig:BE}.
Evidently, the $r_n(t;0)$ are bounded for all $t \geq 0$, which
provides evidence that the approximation~(\ref{eq:pert}) satisfies~(\ref{eq:v}) uniformly to ${\rm O}(\epsilon^2)$.    The implied constant grows
as $\alpha \to 1$, however, as can be seen
by examining the  
vertical scales in Figure~\ref{fig:BE}.
Therefore, the usefulness of~(\ref{eq:pert}) as a 
source of initial guesses to the numerical method
diminishes accordingly as $n$ increases. 

\begin{figure}[t!] 
\centering
\hspace*{-6pt}
    \includegraphics[height=95pt]{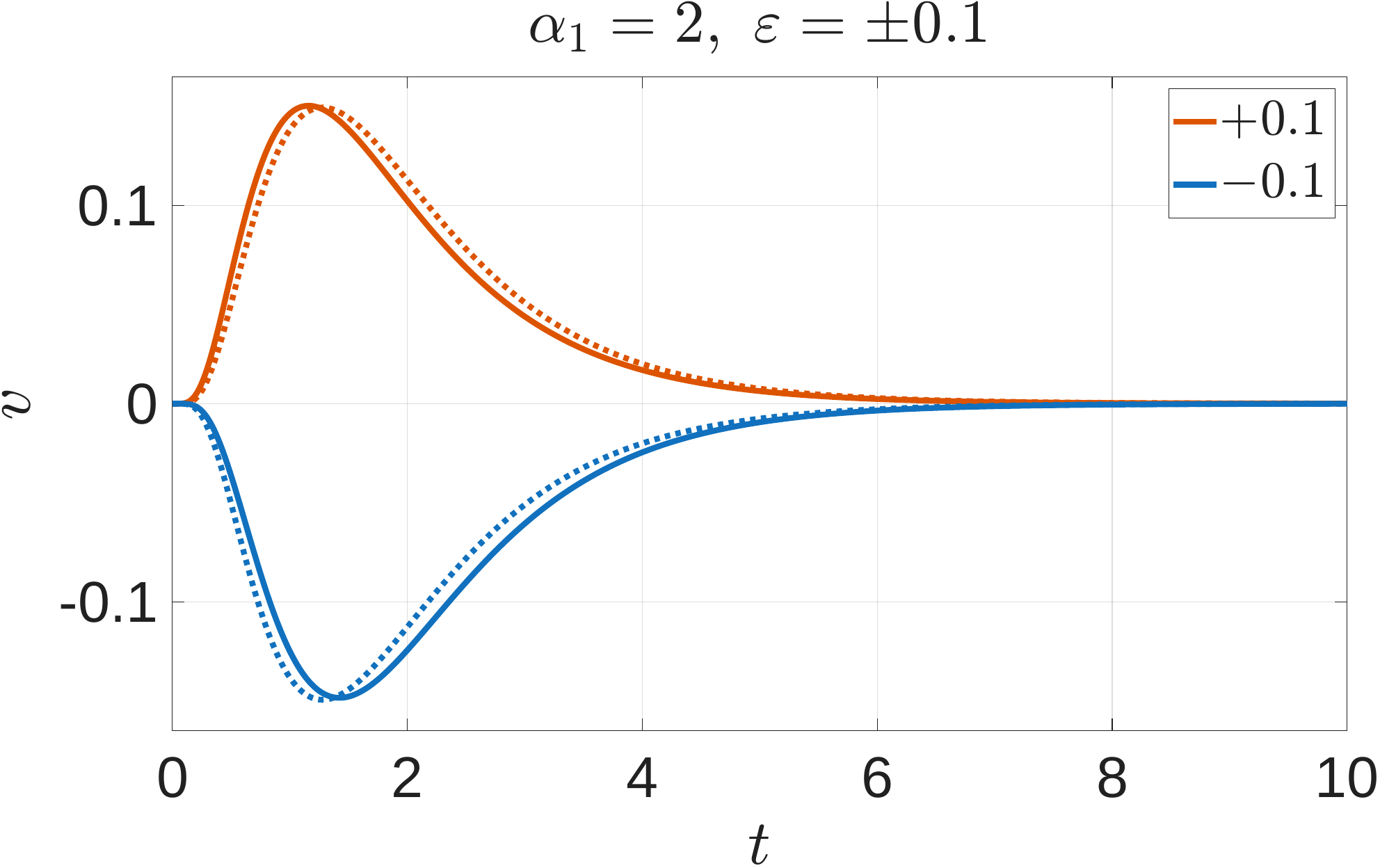}%
    \includegraphics[height=95pt, trim={22pt 0 0 0},clip]{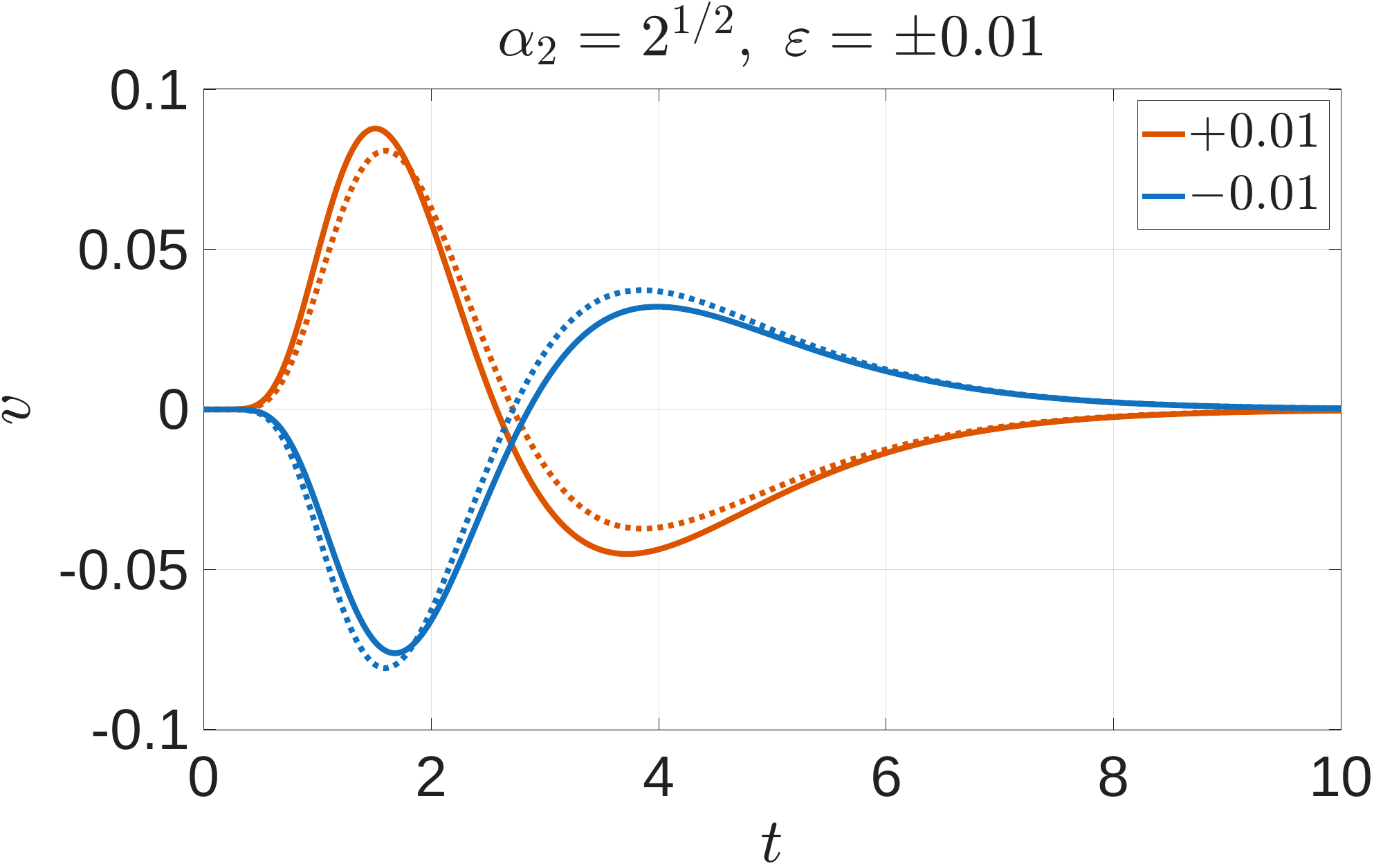}%
    \includegraphics[height=95pt, trim={22pt 0 0 0},clip]{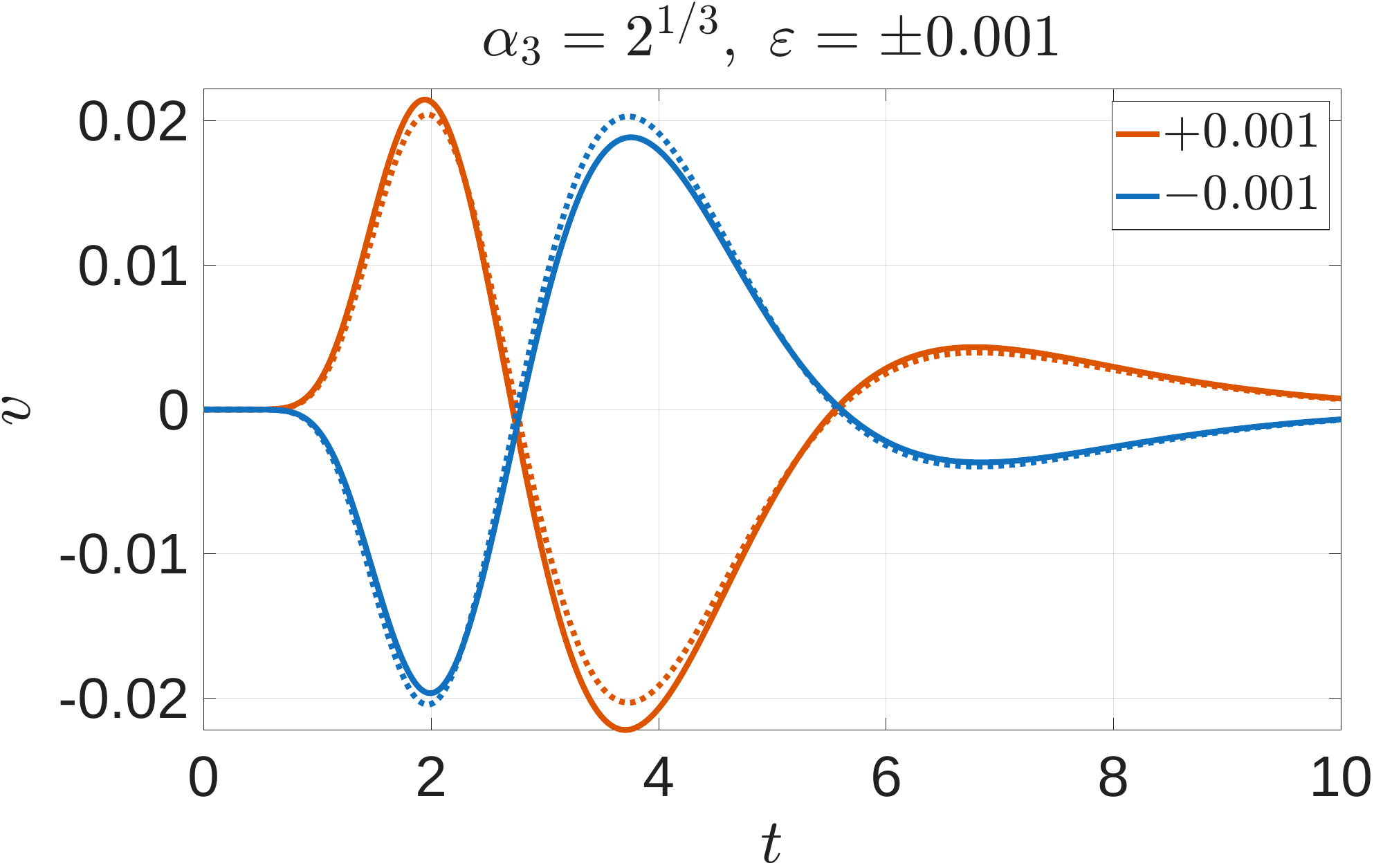}\\[-1em]
    \caption{Numerical solutions $v$ of the perturbation 
    equation~(\ref{eq:v}) (solid curves) in comparison to the
    perturbation approximations~(\ref{eq:pert})~(dashed curves). The scaling constants $C_1$, $C_2$, and $C_3$ are given in Table~\ref{table:c}.
    In each frame the value of $\alpha$ is $\alpha_n+\epsilon$,
    with numerical values of these parameters given in the titles and legends.}
    \label{fig:pertVSnum}
\end{figure}

The details of the numerical method are postponed
to the next section, but let us now 
solve~(\ref{eq:nonlinear}) (via~(\ref{eq:v}))
for values of $\alpha$ close to the
characteristic values. 
A few such solutions
are shown in~Figure~\ref{fig:pertVSnum}, where the
perturbation and numerical approximations 
are displayed.  Evidently, the perturbation
approximations provide good qualitative descriptions of 
the solutions. 

We conclude this section with two remarks.
First, the $v$ solutions shown in~Figure~\ref{fig:pertVSnum} can be distinguished
according to whether the first turning point
is a local maximum (shown in red) versus a~local minimum (shown in blue).   This will continue to be a distinguishing
factor in the solutions presented below; see for example
Figures~\ref{fig:subplot2} and~\ref{fig:subplot3} 
in section~\ref{sec:family}. 
Second, note that the solutions shown 
in~Figure~\ref{fig:pertVSnum} may not be the only
ones for the given parameter values.  (Additional
solutions for $\alpha = \sqrt{2}\pm 0.01$
are shown in subplots {\tt J} and {\tt L}
of Figure~\ref{fig:subplot3}.)
These additional solutions can be produced by starting from the
ones shown in~Figure~\ref{fig:pertVSnum} and then
using a continuation procedure as  
described in the next section.

\section{Numerical method} \label{sec:nummeth}

Motivated by the semi-infinite domain, the smoothness and exponential decay of the solutions, and a desire for high accuracy, we solve~(\ref{eq:v}) numerically using a Laguerre spectral collocation method with scaling appropriate to match the $v(t)={\rm O}(e^{-t})$ behaviour for large $t$ \red{(see~\cite[Chapter 17.5]{boyd}, \cite[Section 2]{dmsuite})} \red{and functional terms implemented via the approach described in~\cite{Hale2023}}. To improve efficiency {and stability} we follow~\cite{Mastroianni} in using 
a truncated Laguerre method. In particular, here we approximate the solution on   
the first $\red{m =} \lceil M\sqrt{N}\rceil$ points of an $N$-point Laguerre grid, with $M$ determined heuristically based on experiments below. 

\red{
Let us make this discretisation more explicit. For a fixed scaling parameter $\beta>0$, the Laguerre spectral approximation is written as a weighted degree-$N$ polynomial, 
\be 
v_N(t)=\exp(-\beta t/2)p_N(t).\ee
In the computations below we take $\beta=2$, matching the expected $v(t)={\rm O}(e^{-t})$ decay. In Laguerre spectral collocation, the polynomial $p_{N}$ is represented in barycentric-Lagrange form using $0 \cup \bm{t}$ as support points, where $\bm{t} = [t_1,  \ldots, t_{N}]^\top$ denotes the zeros of the scaled degree $N$ Laguerre polynomial, $L_{N}(\beta t)$. The point $t=0$ is included to facilitate the implementation of the boundary condition at the origin. The corresponding {differentiation matrix} is determined by exactly differentiating the interpolant $v_N(t)$ and evaluating at the same support points. Then, since here $v(0) = 0$, the first row and column of the differentiation matrix can be removed, resulting in an $N\times N$ matrix,  $D$, such that 
$
v'(\bm{t}) \approx Dv(\bm{t}).
$
The points $\bm{t}$ and the matrix $D$ can be constructed numerically using, for example, {\tt lagdif} in DMSUITE~\cite{dmsuite}, or the more stable implementation in~\cite{Nel2026}.}

The functional terms, $v(\alpha t)$ and $v^2(\alpha t)$, in the equation are implemented with the barycentric resampling approach described in~\cite{Hale2023}. 
\red{
Let $P(\bm{s},\bm{t})$ denote the barycentric interpolation matrix that maps values of an unweighted polynomial at the nodes $\bm{t}$ to values at the target points $\bm{s}$.  The resampling matrix used to evaluate the weighted interpolant at the dilated nodes $\alpha\bm{t}$ is
\be
P_\alpha =
\operatorname{diag}\!\left(e^{-\beta\alpha \bm{t}/2}\right)
P(\alpha\bm{t},\bm{t})
\operatorname{diag}\!\left(e^{\beta\bm{t}/2}\right).
\ee
Thus, if  $\bm{v}\approx v(\bm{t})$ is the vector of approximate solution values on the Laguerre grid, then $P_\alpha\bm{v}$ approximates $v(\alpha\bm{t})$. 
The matrix $P_\alpha$ can be constructed numerically using {\tt barymat} in Chebfun~\cite{Chebfun}.
}

Hence, the full Laguerre spectral discretisation of~(\ref{eq:v}) yields the
\red{$N\times N$} nonlinear system
\be\label{eq:discretisation} (D + I - 2P_\alpha)\bm{v} - \red{(P_\alpha\bm{v})\circ(P_\alpha\bm{v})} = 0,\ee
where $\circ$ is the element-wise product.

\begin{figure}[t!]  
\centering
    \includegraphics[height=110pt]{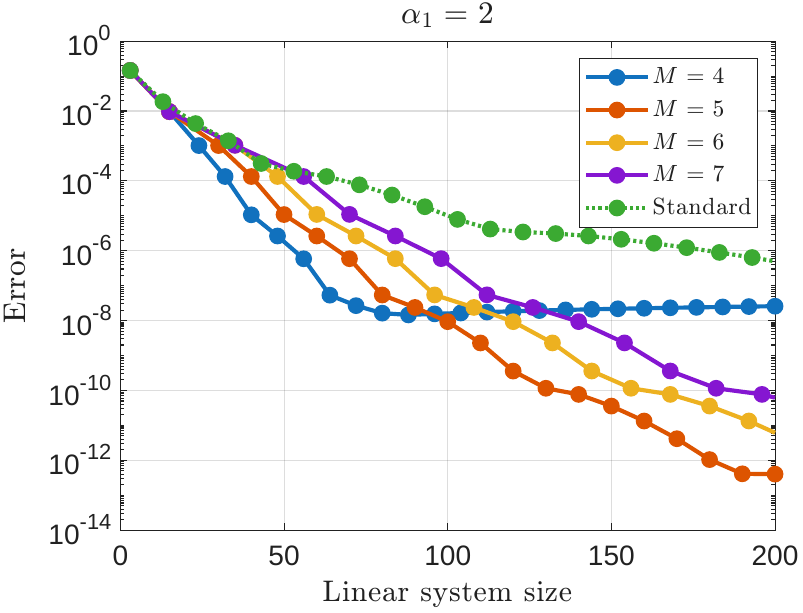} \ %
    \includegraphics[height=110.75pt, trim={0 0 0 0},clip]{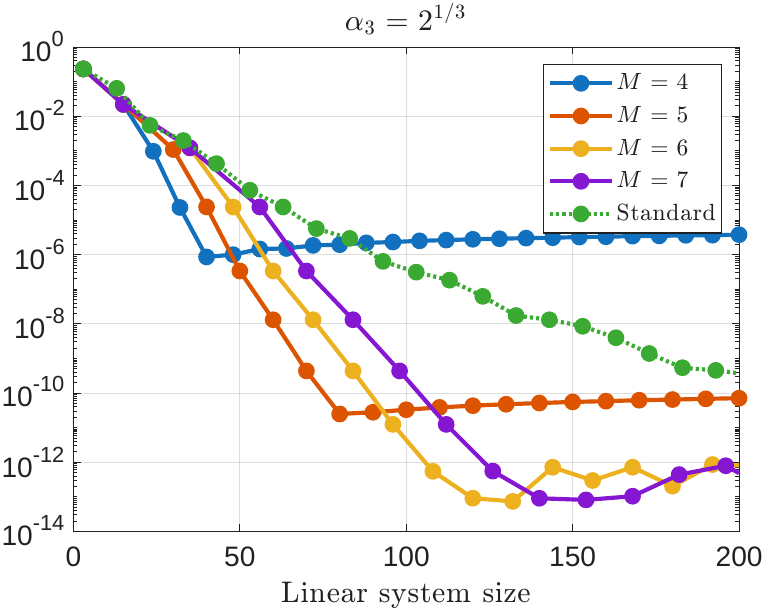} \ %
    \includegraphics[height=110pt, trim={0 0 0 0},clip]{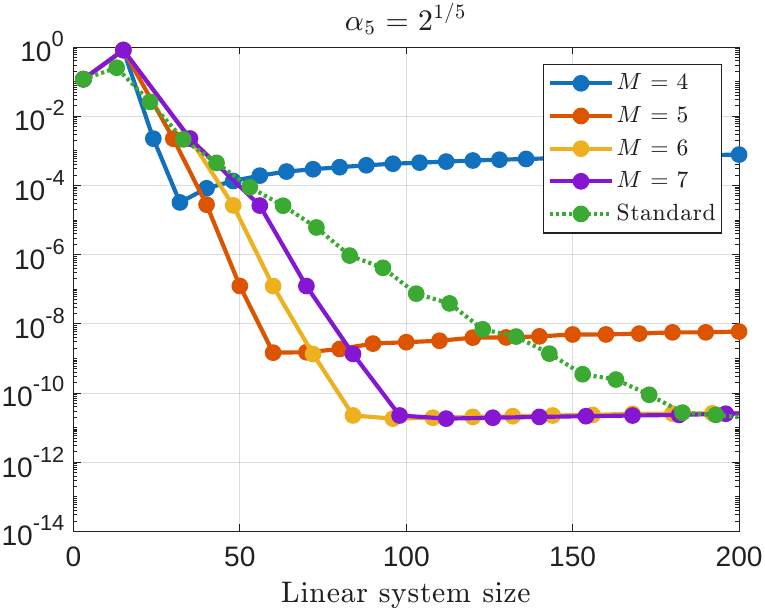}\\[-1em]
    \caption{Errors in the $\red{m = }\lceil M\sqrt{N}\rceil$-truncated Laguerre spectral discretisation of~(\ref{eq:linear2}) for the characteristic solutions {as a function of discretisation size}, for $\alpha = 2, 2^{1/3},$ and $2^{1/5}$. Shown  are the  $\ell_2$
    grid errors as compared to the series solution~(\ref{eq:series}),  \red{normalised
     as in Figure~\ref{fig:char}}. {The curve marked `Standard' depicts the same errors when using the standard Laguerre discretisation (i.e., no truncation). The truncated approach with  $M=6$ provides the best balance between convergence and stability.}
    \red{The final panel indicates that the attainable accuracy begins to be limited by the conditioning of the discretised system as $\alpha\to1$.}
    }
    \label{fig:verifylinear}
\end{figure}

\red{The truncated version is similar, but retains only the first $m=\lceil M\sqrt{N}\rceil$ nodes of $\bm{t}$. One can achieve this by constructing  $D$ and $P_\alpha$ in full and retaining only their first $m\times m$ sub-blocks, or by directly computing the truncated differentiation and barycentric resampling matrices corresponding to 
a weighted degree $m$ interpolant 
with support points $[0,t_1,  \ldots, t_{m}]$. 
We choose the former for simplicity. In either case, one again arrives at~(\ref{eq:discretisation}), but now the truncated system is only $m\times m$.} {Figure~\ref{fig:verifylinear} shows that this approach improves convergence (in terms of error versus the size of the linear/nonlinear systems to be solved) compared to \red{the standard approach of} using all $N$ points.}

We begin by verifying the efficacy of the Laguerre discretisation by solving the linear problem~(\ref{eq:linear2}) and comparing to the series solution~(\ref{eq:series}). To do so, we omit the nonlinear term from~(\ref{eq:discretisation}) and obtain nontrivial solutions by computing the 
nullspace of the matrix $D + I - 2P_\alpha$ \red{for $\alpha = \alpha_n$}. Figure~\ref{fig:verifylinear} shows the error between the spectral collocation and series solutions for various $\alpha_n$ as the 
number of degrees of freedom
in the collocation method is increased. 
\red{
In particular, the error plotted in Figure~\ref{fig:verifylinear} is the grid error
\be
    \|E(\bm{t};\alpha_n)-\gamma\bm{v}\|_2,
    \qquad
    \gamma=\frac{E(\bm{t};\alpha_n)^\top\bm{v}}{\bm{v}^\top\bm{v}},
\ee
where $E$ is normalised as in Figure~\ref{fig:char}. The scalar $\gamma$ is the least-squares scaling factor for the computed null vector $\bm{v}$, whose normalisation is arbitrary. The reference values of $E(\bm{t};\alpha_n)$ in~(\ref{eq:series}) are computed directly at the retained Laguerre nodes by summing the series in extended precision. Extended precision is necessary here because, although the terms eventually decay rapidly, the coefficients may undergo substantial transient growth and cancellation, especially for larger $n$, resulting in inaccurate computation in double precision arithmetic. 
}

The error curves indicate rapid convergence up to a point where \red{either} a sub-optimal truncation point degrades the quality of the solution \red{or the attainable accuracy is limited by the deteriorating conditioning of the discretised linear system}.
{The dotted line depicts the standard Laguerre discretisation (no truncation) and we see that the truncated versions are more efficient.}
Rather than trying to find the optimal scaling and discretisation size for each $\alpha$, based on the results in Figure~\ref{fig:verifylinear}, we take $M=6$ and $N = 700$ (so that $\red{m} = 159$)  in our remaining experiments to provide a balance between accuracy and execution time.
\red{
For these final parameter values, the $\ell_2$ grid errors for the linear problem are reported in Table~\ref{tab:linear-errors}.
}

\begin{table}[h!]
\begin{center}
\red{
\begin{tabular}{c|cccccc}
$\alpha$ 
& $2$ 
& $2^{1/2}$ 
& $2^{1/3}$ 
& $2^{1/4}$ 
& $2^{1/5}$ 
& $2^{1/6}$ \\ \hline
Error
& $1.1\times10^{-10}$ 
& $3.4\times10^{-14}$ 
& $1.7\times10^{-13}$ 
& $4.7\times10^{-13}$ 
& $2.4\times10^{-11}$ 
& $4.5\times10^{-10}$
\end{tabular}
}
\end{center}
\caption{\red{Grid errors, $\|E(\bm{t},\alpha)-\gamma\bm{v}\|_2$, for the linear problem~(\ref{eq:linear2}) using $M=6$ and $N=700$ ($m=159$) in the truncated Laguerre spectral collocation method. These parameter values are used to solve the full nonlinear problem~(\ref{eq:nonlinear}) in the next section.}}
\label{tab:linear-errors}
\end{table}

Having verified that the Laguerre discretisation 
provides accurate solutions to the linear problem~(\ref{eq:linear2}), we now use it to solve the nonlinear problem~(\ref{eq:v}). We first do so in the vicinity of the characteristic values, $\alpha_n=2^{1/n}$. In particular, we use the perturbation approximations derived in section~\ref{sec:nonlinear} as initial guesses 
which are then refined using a Newton 
iteration applied
to~(\ref{eq:discretisation}). Some examples have already
been shown in Figure~\ref{fig:pertVSnum}. For each $n$, these solutions are then used to seed a pseudo-arclength continuation, 
as described by \cite{Farrell2015} and implemented in Chebfun~\cite{Chebfun},
to compute solutions over a desired range of $\alpha$. 
 
A MATLAB code that implements the method described here
is available at~\cite{code}.  It can be used to reproduce 
all figures of this paper.

\section{A family of solutions} \label{sec:family}

Using the numerical procedure described in the previous section,
we now  explore the family of solutions in class (A) as defined in the first paragraph
of section~\ref{sec:nonlinear}, i.e.,  the solutions to~(\ref{eq:nonlinear}) that satisfy $u(0) = 1$ and $u(t) = 1 + {\rm O}(e^{-t})$ as $t \to \infty$.  Our aim is to catalogue for each value of $\alpha$ the number of possible 
solutions in this class.  Recall that 
for each value of $\alpha$ the constant solution $u = 1$ is in this class, so we focus on counting the additional solutions. 
Excluded in the counting is the solution that satisfies $u(t) = {\rm O}(e^{-t})$ as $t \to \infty$,
as it is not in class (A).  Possible solutions that obey an asymptotic behaviour different
from ${\rm O}(e^{-t})$  are excluded in the counting as well.

\begin{figure}[thbp]   
\begin{center}
    \hspace*{-1pt}\includegraphics[width=.452\textwidth]{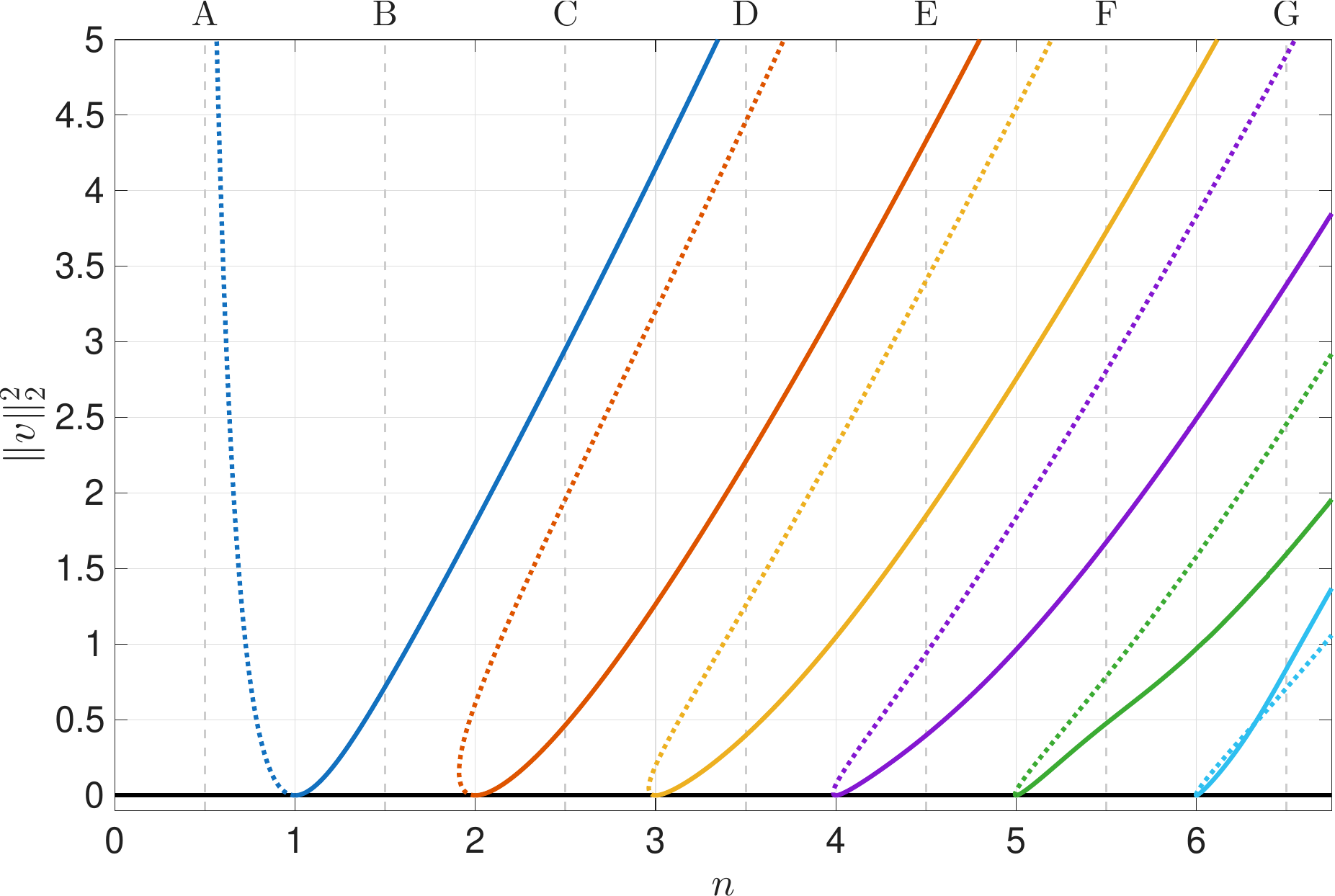}\qquad\hspace*{-4pt} \includegraphics[width=.45\textwidth]{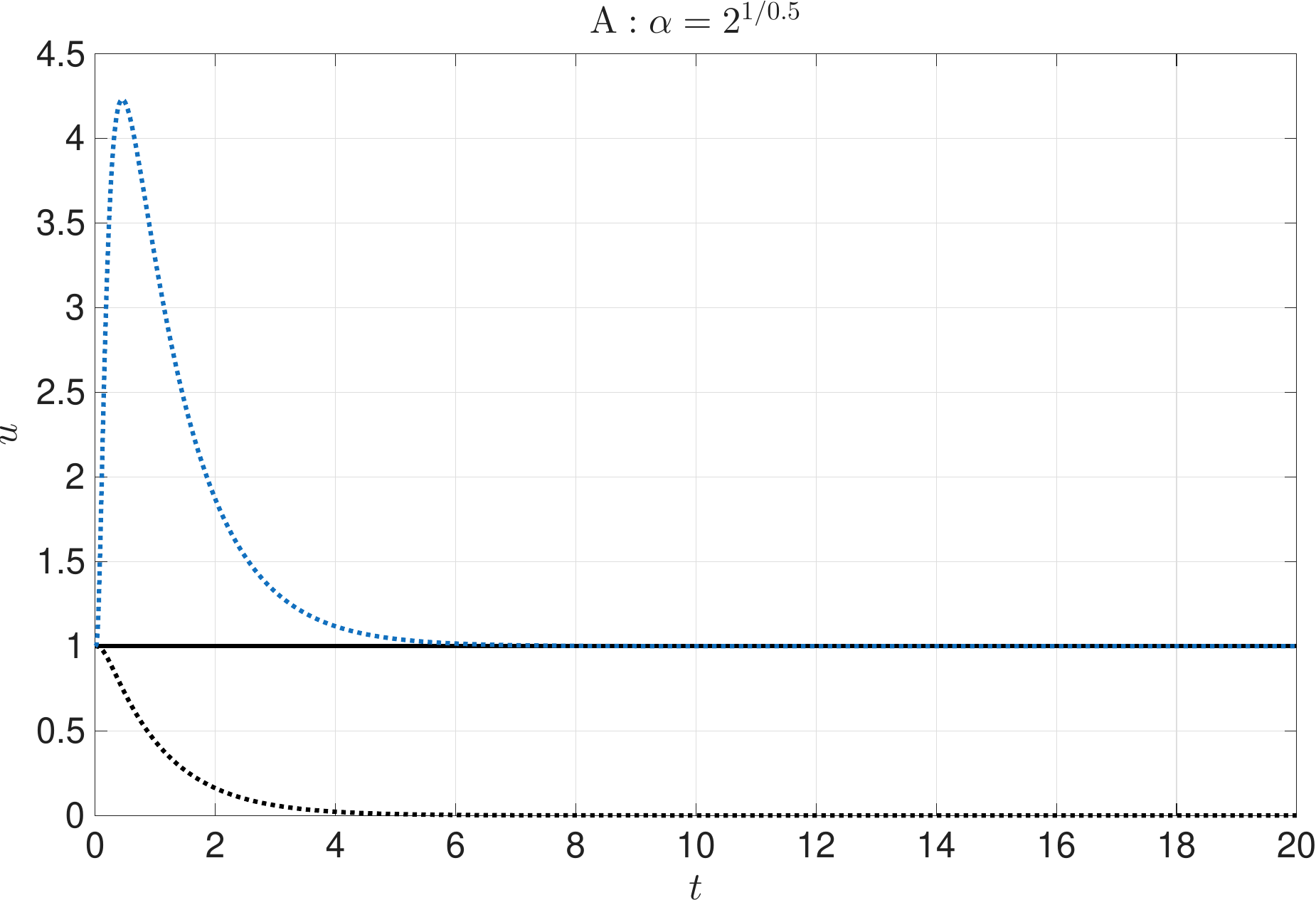}\\    \includegraphics[width=.45\textwidth]{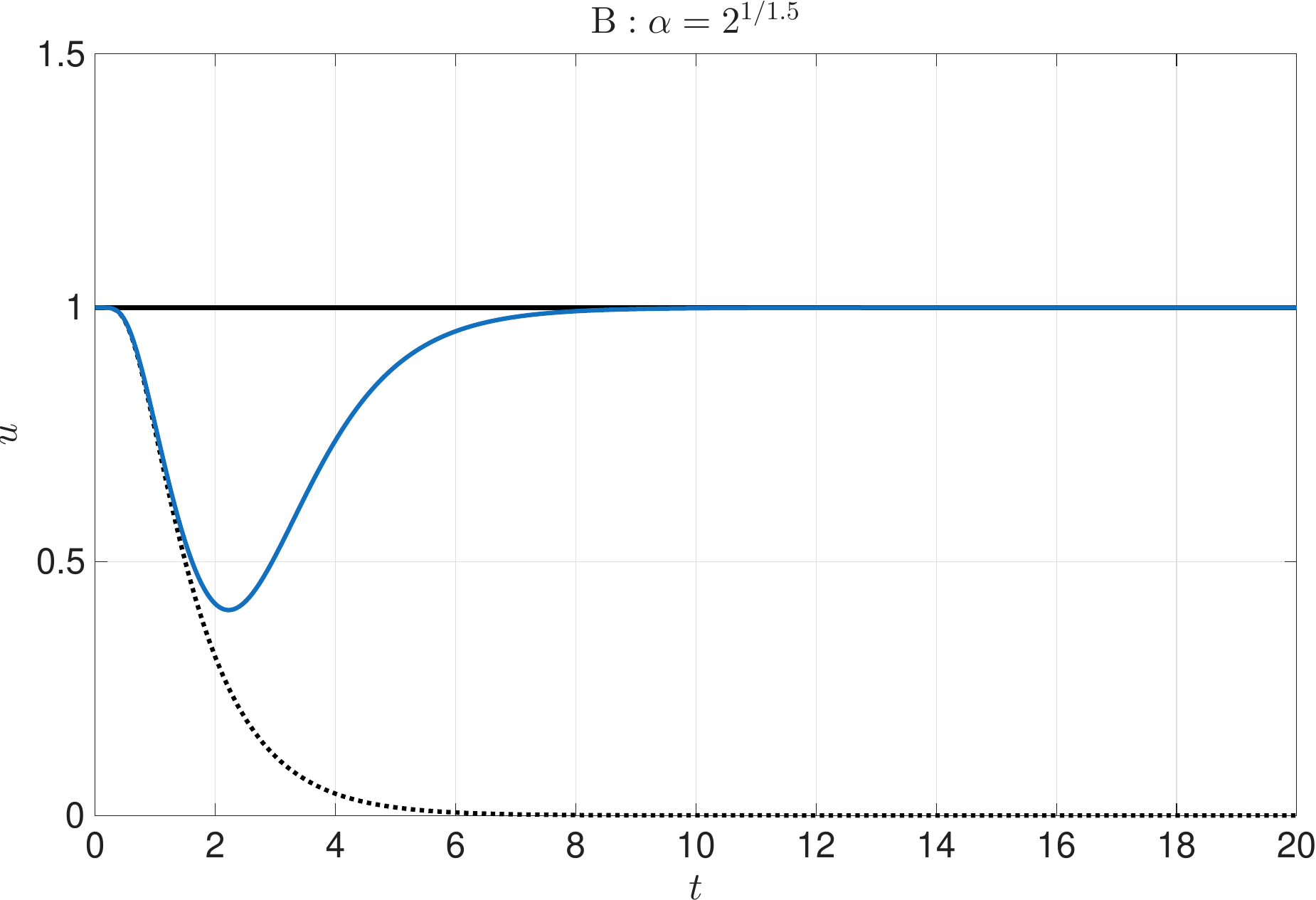}\qquad
    \includegraphics[width=.45\textwidth]{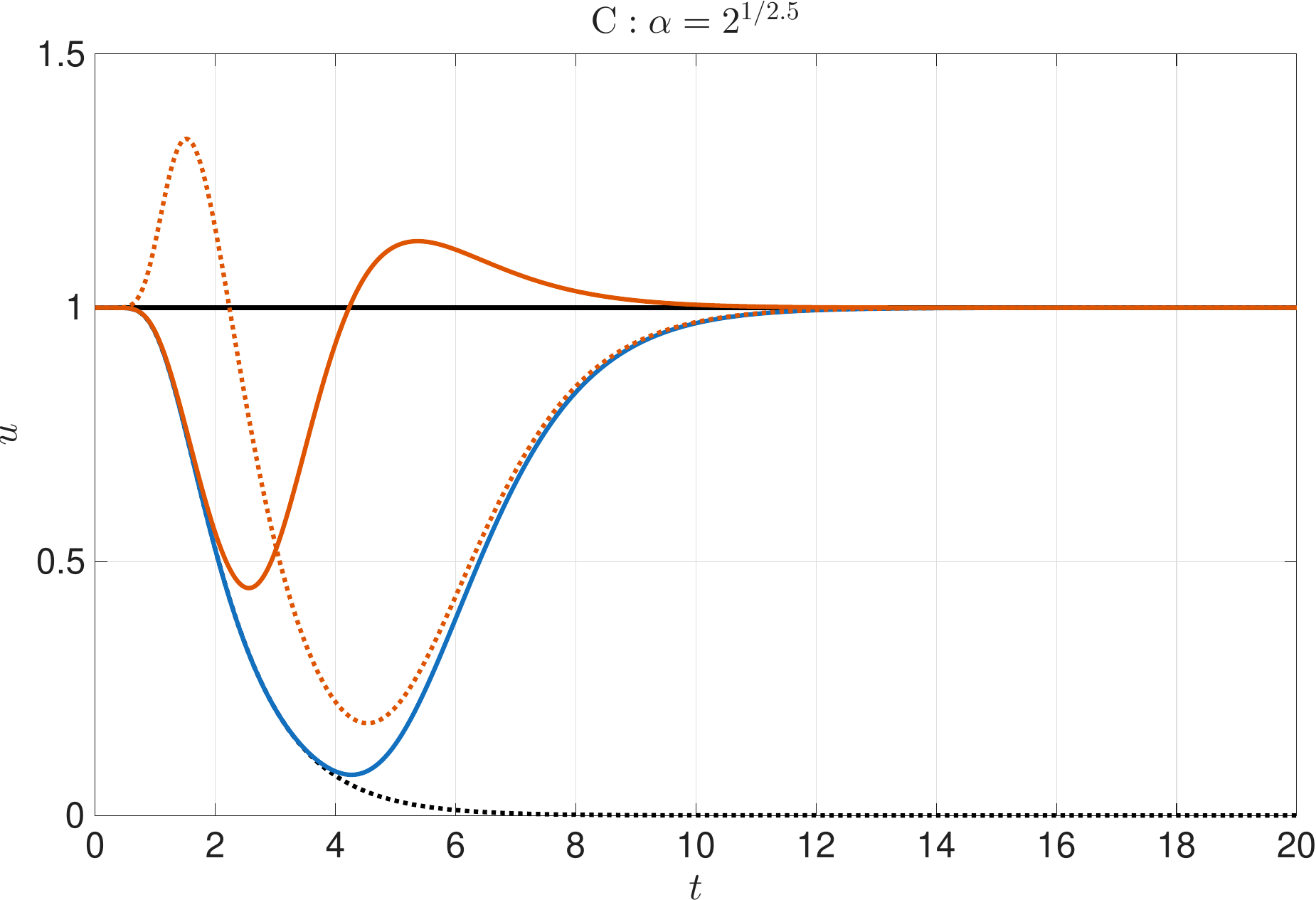}\\
    \includegraphics[width=.45\textwidth]{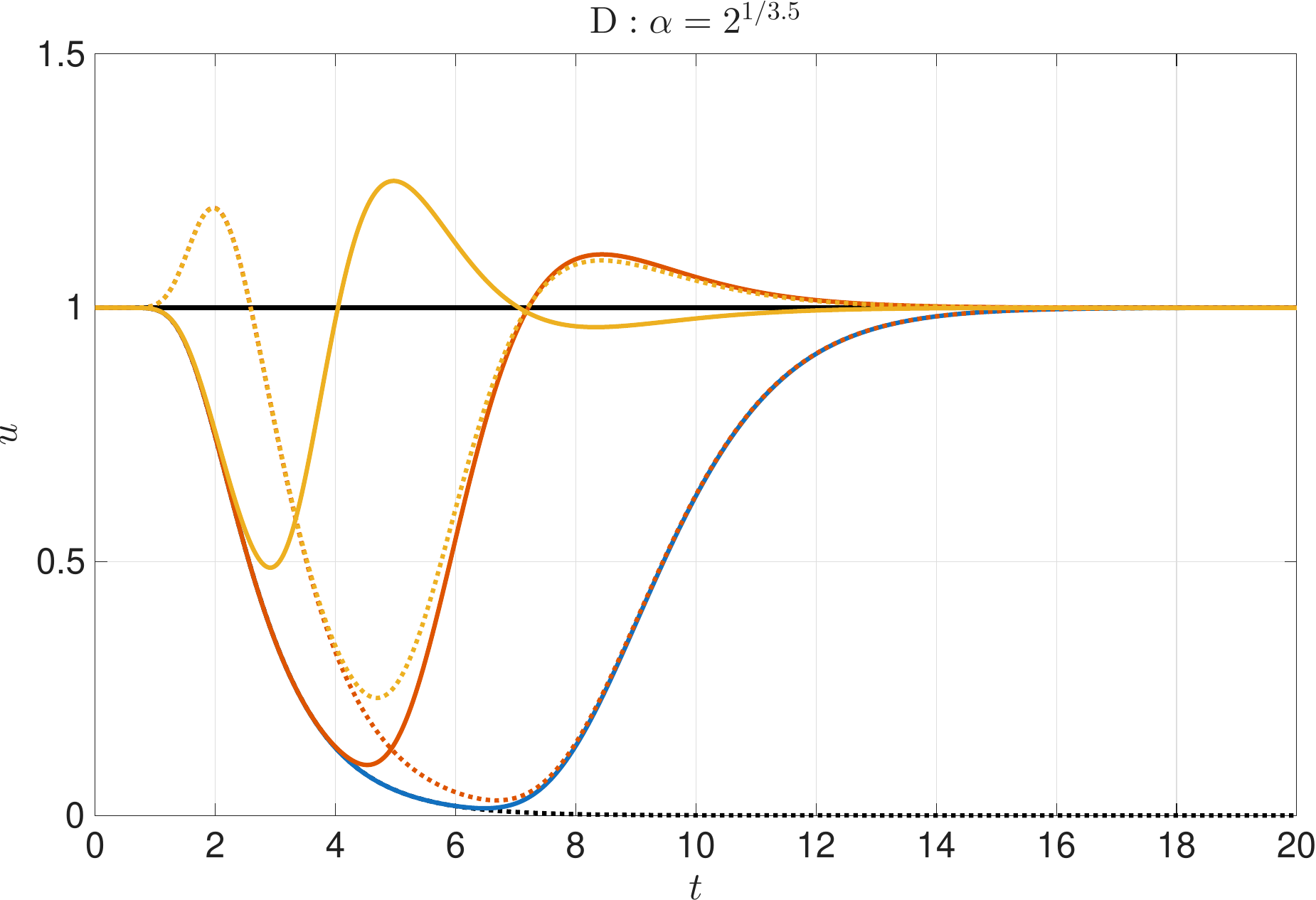}\qquad
    \includegraphics[width=.45\textwidth]{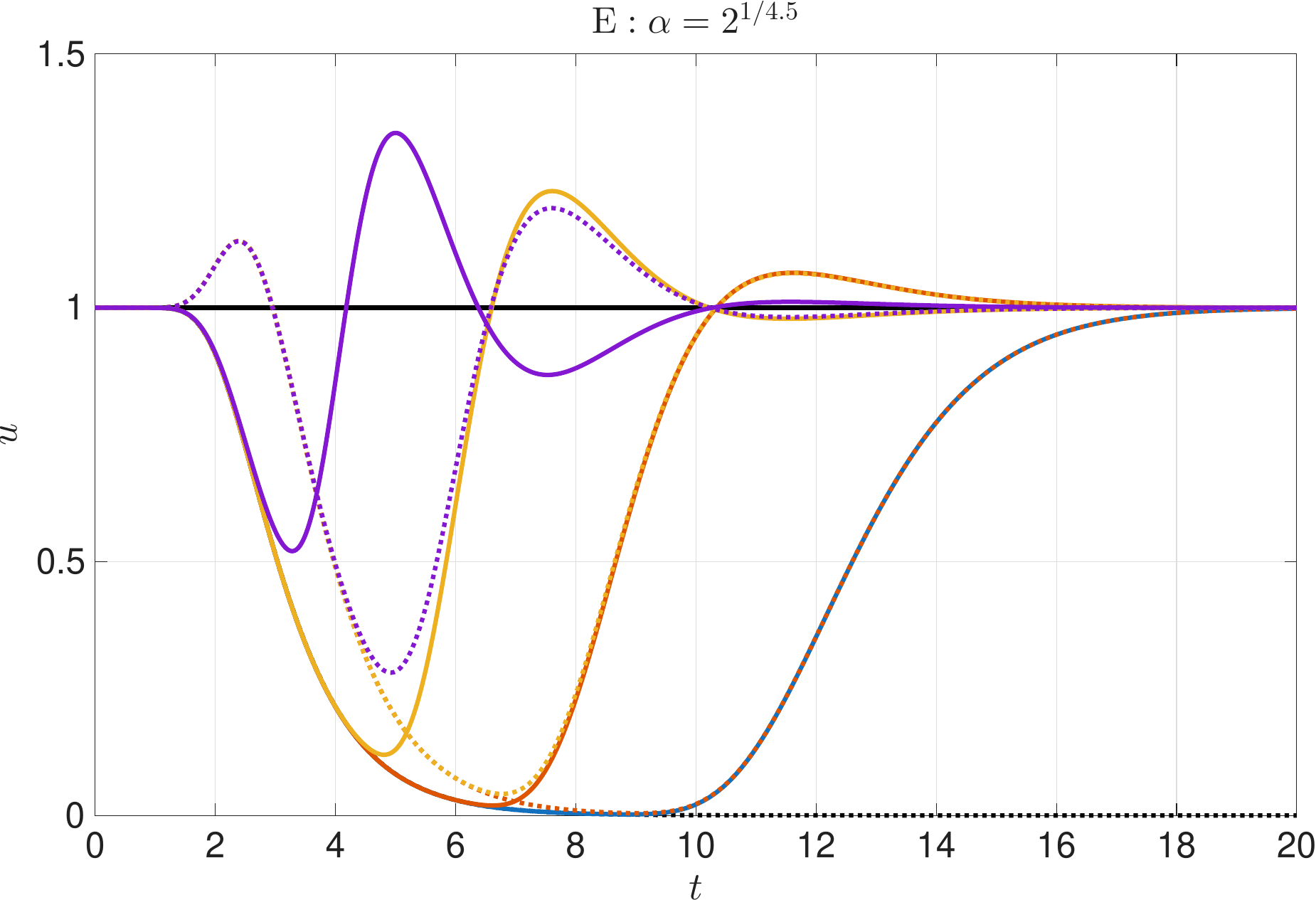}\\
    \includegraphics[width=.45\textwidth]{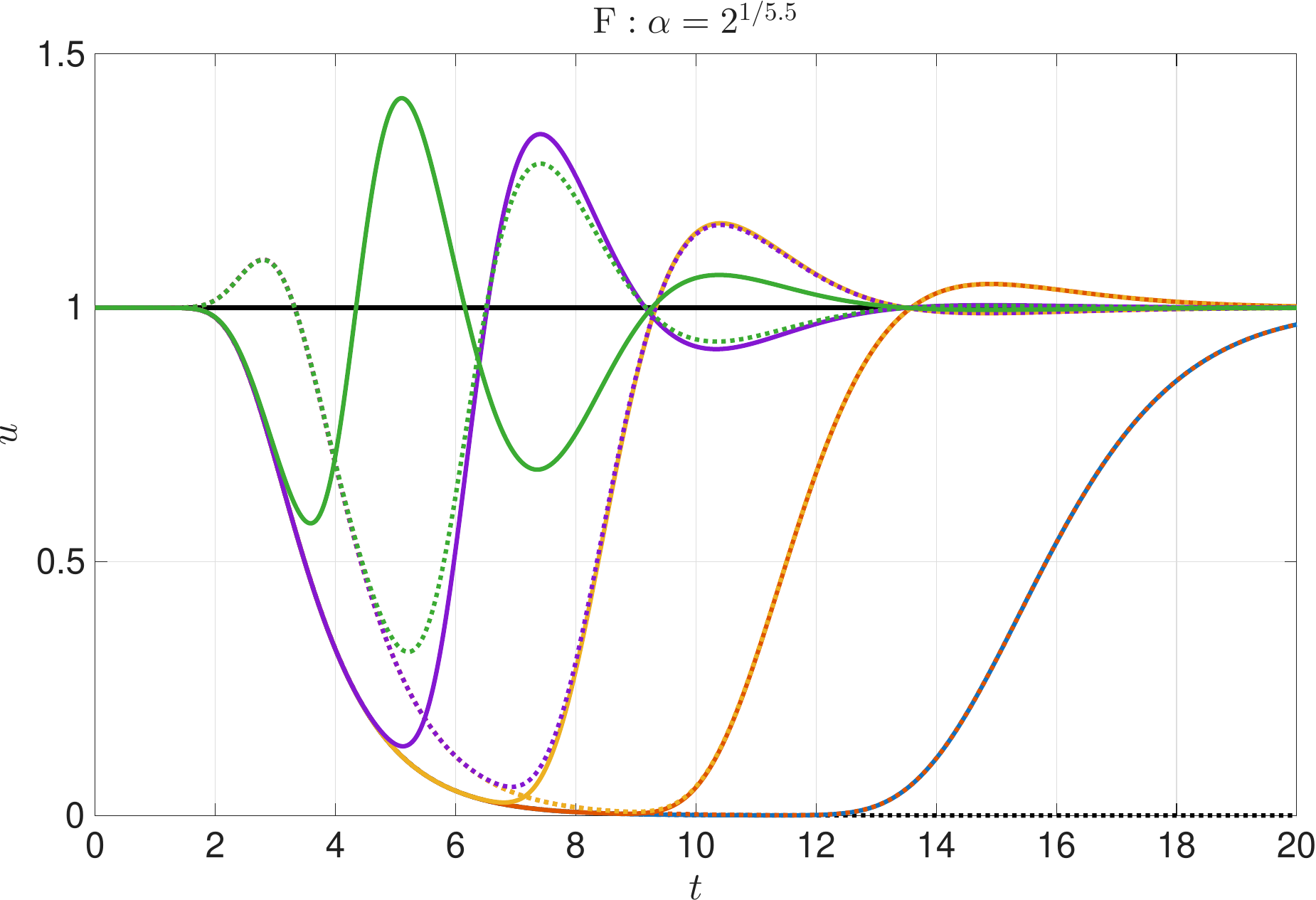}\qquad
    \includegraphics[width=.45\textwidth]{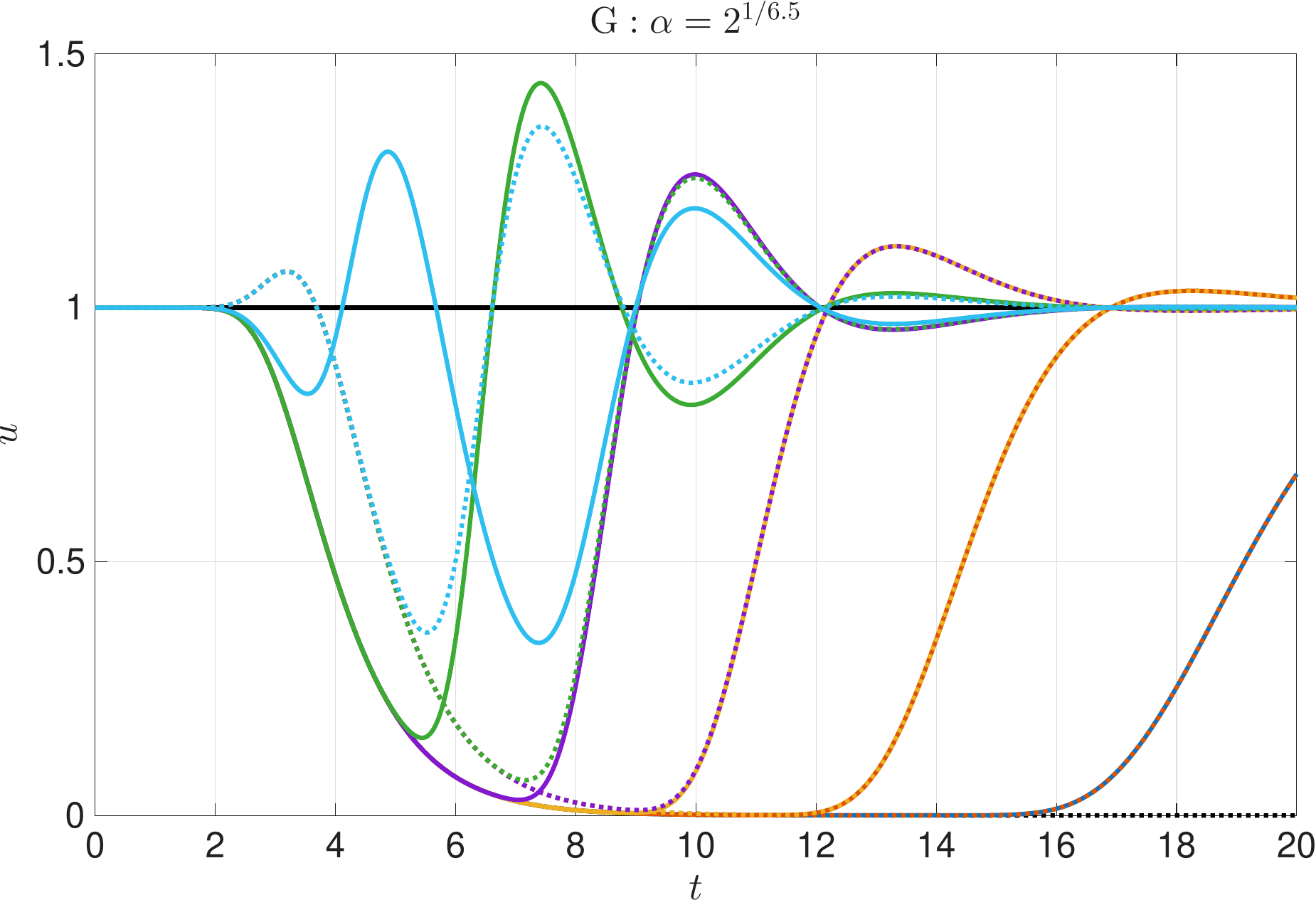}
\end{center}

\vspace*{-1em}
    \caption{The first frame shows the bifurcation diagram and the remaining frames show
    solutions to~(\ref{eq:nonlinear}) for various values of $\alpha$. The black dotted
    solution decaying to zero 
    is the one in class (B) as defined in the 
    first paragraph of section~\ref{sec:nonlinear}.  The
    multicoloured solutions are in class (A).  Note the labeling:
the letters at the top of the bifurcation diagram correspond to the various
solution plots.     Dotted (resp.~solid) curves are used to 
represent solutions whose first turning point is a local
maximum (resp.~minimum). 
    (An animated version of the figure can be seen
    at \href{https://youtu.be/aBTvLujm65c}{https://youtu.be/aBTvLujm65c}.)}
    \label{fig:subplot2}
\end{figure}

A bifurcation diagram for solutions in this family is shown in the first frame of Figure~\ref{fig:subplot2}, and solutions to~(\ref{eq:nonlinear}) are shown in the other frames. \red{The solutions in this figure are displayed on the interval $[0,20]$ only for visualisation; this plotting window is not a finite-domain truncation in the Laguerre discretisation, which is posed on $[0,\infty)$.
}
The bifurcation diagram shows the values of the   $L^2$-norm (squared) of solutions 
$v$ to~(\ref{eq:v}) as $\alpha$ is varied.   
Points at which the curves touch the horizontal axis correspond to $v=0$ and therefore represent the constant solution $u = 1$.  This happens at the 
 critical values of $\alpha$, namely $\alpha_n = 2^{1/n}$,
 $n = 1, 2, \ldots $.  In the vicinity of these
 points, the continuation procedure described in the previous section was used to generate further solutions.

 Consider first the case $n = 1$ (i.e., $\alpha_1 = 2$). 
 Precisely at this value, the only solution
 in class (A) is the constant solution
 $u=1$, but non-constant solutions emerge if $\alpha$
 is perturbed slightly. Such solutions, corresponding
 to $\alpha = 1.9$ and $2.1$, were already displayed in the
 first frame of Figure~\ref{fig:pertVSnum}  
 (under the translation $u=v+1$).   At the end of 
 section~\ref{sec:nonlinear} we made a remark about
 distinguishing such solutions according to whether they
 first dip down or flip up as $t$ is increased. This
 leads us to the following conventions: Solutions $u$  whose first turning point on the positive $t$-axis 
 is a local maximum will be called {\em plus-solutions}
 and their graphs will be plotted with dotted curves. 
 If the first turning point is a local minimum, the solutions are {\em minus-solutions} and are plotted with solid curves.

 Let us now consider {values of $\alpha$ away from the critical points, corresponding to noninteger $n = \log(2)/\log(\alpha)$.} First, suppose one decreases
 the value of $n$ from $1$ down toward $0$ (i.e., $\alpha$ increases from the value $2$).
 Following the first bifurcation curve to the left (dotted blue), one encounters for each value of $\alpha$ a unique 
 plus-solution.  An example is the {$n = \frac12$} ($\alpha = 4$) solution
 shown in the
 second frame of~Figure~\ref{fig:subplot2}.      Increasing
 the value of $\alpha$ further  leads to similar plus-solutions with increasing norm. 

 Next, consider increasing the value of $n$ from $1$ (i.e., $\alpha$ decreases from $2$). 
 Following the first bifurcation curve to the right 
 (solid blue), one encounters for each value of $\alpha$
 a unique minus-solution.  An example can be seen in the third frame  of~Figure~\ref{fig:subplot2}.
 When $n$ gets sufficiently close to $2$ ($\alpha$ close to $\sqrt{2}$),  however, 
 two additional solutions emerge as represented by the red curve in the bifurcation diagram. 
 Now there are three solutions:
 the minus-solution of large norm represented by the blue curve as well as the two new solutions, 
 both with relatively small norm.   

\begin{figure}[t]   
\centering
    \hspace*{-3pt}\includegraphics[width=.324\textwidth,trim={1pt -1pt 2pt 0pt},clip]{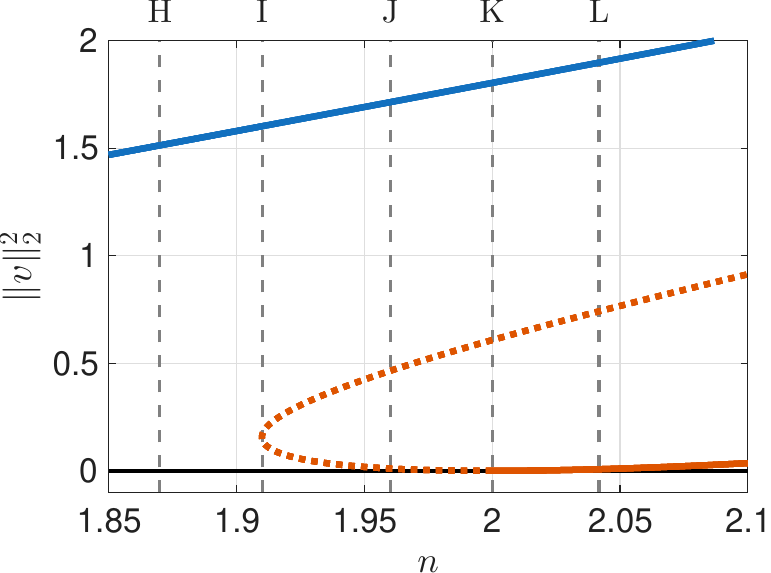}
    \includegraphics[width=.32\textwidth]{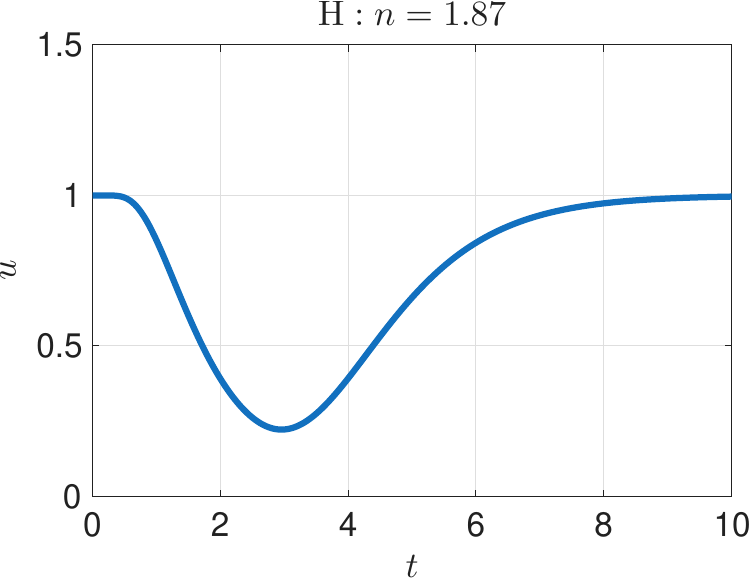}
    \includegraphics[width=.32\textwidth]{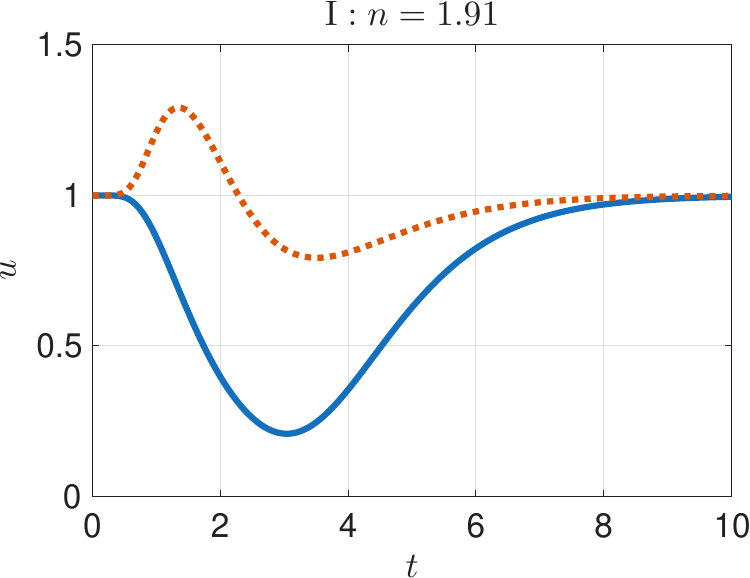}\\
    \includegraphics[width=.32\textwidth]{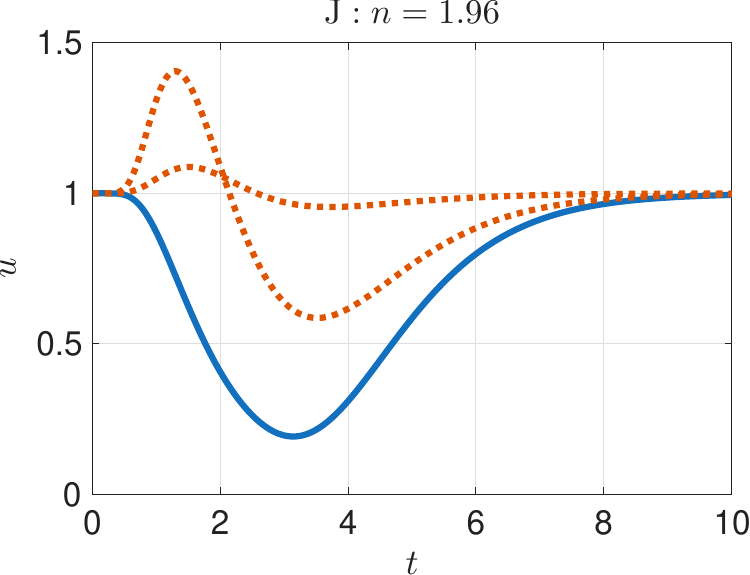}
    \includegraphics[width=.32\textwidth]{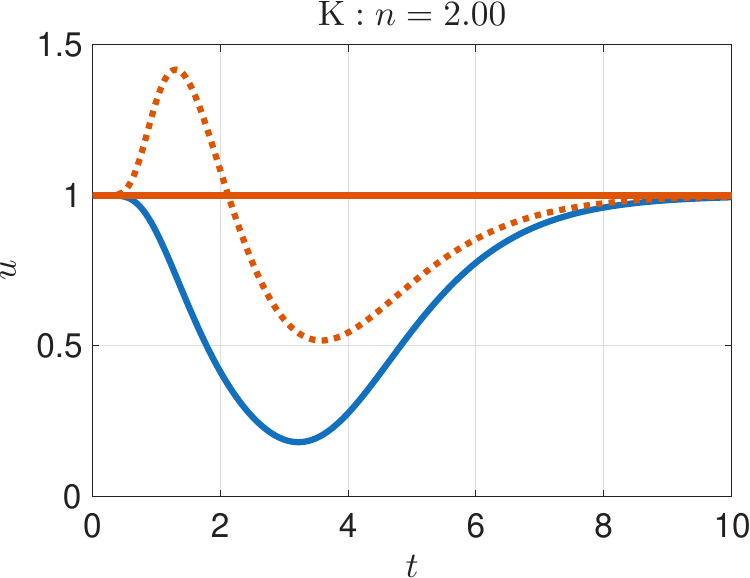}
    \includegraphics[width=.32\textwidth]{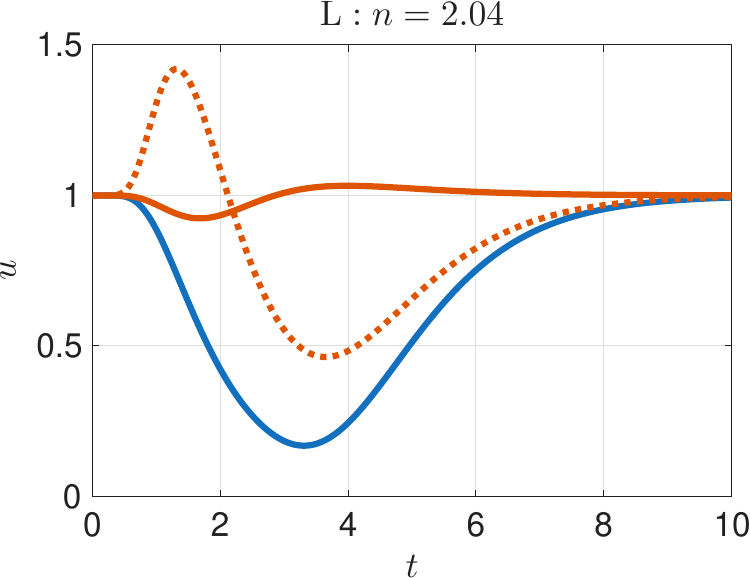}\\[-1em]
    \caption{The bifurcation diagram on the top left is a small section of the
    one shown in the first frame of Figure~\ref{fig:subplot2}.   The  other frames   show the generation of new solutions
    as the value of $n$ is increased through the value $2$ (i.e., as
    $\alpha$ is decreased through the value~$\sqrt{2}$).}
    \label{fig:subplot3}
\end{figure}

For more detail about these solutions, 
consider Figure~\ref{fig:subplot3}, which zooms in 
on the interval around $n =2$ in the bifurcation diagram.
When $n$ is less than $2$ by a sufficient margin, 
the only solution in class~(A) (other than the 
constant solution) is a minus-solution such as the one
shown with the label~{\tt H}.   (This solution is
similar to the one computed  
in~\cite[Figure~2]{DascaliucErratum} for a slightly different value of $\alpha$.)
The solution
marked {\tt I} corresponds to the critical value of $n \approx 1.91$  ($\alpha \approx 1.44$)
where a first additional solution is generated.   With increasing
$n$ this solution bifurcates into two plus-solutions as shown
in the subplot labeled {\tt J}.   When $n$ reaches the value $2$, 
one of these plus-solutions (the one with the smaller norm) 
momentarily collapses to the constant solution at the bifurcation point (subplot {\tt K}) and
then re-emerges as a minus-solution for larger $n$  
(subplot {\tt L}).  {Our numerical tests suggest that this is indeed a transcritical bifurcation.
	} Observe that the two  solutions of small norm
shown in subplots {\tt J}
and {\tt L} in Figure~\ref{fig:subplot3} correspond
to the two solutions shown in the middle frame of 
Figure~\ref{fig:pertVSnum}.  These are therefore examples of   
the additional solutions alluded to at the end of 
section~\ref{sec:nonlinear}.

Returning to Figure~\ref{fig:subplot2} and 
continuing to increase the value of $n$, one notices a similar pattern. 
In each interval  $(n-\delta, n)$, with $2 \leq n \leq 6$ and 
$0 < \delta \ll 1$, two new solutions are  generated.    
When created, both of these are plus-solutions but one
of them switches to a minus-solution once the value of $n$ is 
exceeded. The magnitude of $\delta$ is about $0.09$ when $n=2$ but then 
decreases rapidly with $n$. 

A noteworthy feature of the family of solutions shown
in~Figure~\ref{fig:subplot2} is the fact that
the minus-solutions all closely follow the exponentially decaying solution (black dots) for
a length of time before turning upwards to approach the value $1$ asymptotically.   
Other solutions have a tendency to merge as well. 
The plus-solution that emerges at integer value $n$ catches up to the minus-solution 
that emerged at integer value $n-1$ as they  jointly approach the limiting value $1$. 

{We conjecture that
the inequality $2 \leq n \leq 6$ above can be replaced by $2 \leq n < \infty$ but numerical experimentation in the range
$n \geq 6.5$ is complicated by a number of issues. First,  the perturbation
approximation~(\ref{eq:pert}) becomes less effective as a starting guess due to the increasing magnitude of the constants $C_n$. Second, increasing oscillations and singular behaviour of the solutions near the origin pose a challenge to the Laguerre spectral discretisation. Third, numerical evidence suggests that the Jacobians in the continuation method become increasingly near-singular, raising questions about the stability of the solutions in this regime.}
 
We conclude with the following puzzling observation. 
As $\alpha \to 1$, more and more solutions {(of increasing $L^2$ norm)} appear.
And yet, when $\alpha$ equals~$1$ exactly, this multitude of solutions collapses
to the unique solution $u = 1$, {suggesting some kind of singular limit}.   At present 
we 
\red{do not have a satisfactory explanation for this} 
nor any means of exploring
further because of the difficulty in computing accurate solutions
in this {limit}, {so the problem remains open}.\footnote{
{The limit $\alpha\to1^+$ in the linear problem was considered in~\cite[Section 3(b)]{Fox} and~\cite[Section 4.4]{Hall}, but it is not clear how those analyses extend to the nonlinear case.}}

\section{Conclusions} \label{sec:conclusions}

The methodology presented here for solving~(\ref{eq:nonlinear}) 
has several ingredients: a Laguerre spectral collocation method for solving the equation numerically,
a perturbation method for generating initial guesses for a  
Newton iteration, and a continuation method for generating further solutions.   
These techniques are well-established for standard ordinary differential
equations, but we believe their combination to solve a nonlinear functional
differential equation is new.

The procedure proved to be a reliable way for generating 
multiple solutions for parameter values $u(0) = 1$ and $\alpha > 1$,
but with $\alpha$ not too close to $1$.  Values of $\alpha$ roughly in
the range $(1,1.11)$ still provide a challenge. {In particular, the initial guesses based on~(\ref{eq:pert}) become less accurate as $\alpha\to1$ ($n\to\infty$) and at the same time the spectral method is adversely affected by both the increasingly oscillatory solutions and the singular behaviour near $t=0$. It is possible that replacing the standard Laguerre basis with generalized Laguerre may lead to improvement, but we have not investigated this.}

\red{
Nevertheless, the approach developed in this manuscript provides a tool for more quantitative investigations of solutions to~(\ref{eq:nonlinear}). For example, Figure~\ref{fig:subplot2} suggests that the length of the initial interval on which some solutions remain close to $1$ increases as $\alpha\to1^+$. For the branch satisfying $0\leq u(t)\leq1$, this behaviour has a probabilistic interpretation, as discussed in~\cite{DascaliucErratum}. Similarly, as $\alpha\to1^+$, solutions appear to remain close to the exponentially decaying class~(B) solution for longer and longer times before eventually returning to $u=1$ as $t\to\infty$. Another observation is the apparent pairing of different solution branches, both for small and large values of $t$. Exploring and understanding these behaviours would be interesting avenues of future work.
}

The methodology should also enable one to consider the following
generalization of~(\ref{eq:nonlinear}),
\be
u^{\prime}(t) + a u(t) = b u^m (\alpha t),
\ee
where $a$ and $b$ are positive constants and $m \geq 2$ is
an integer.    {linearisation about the constant solution
$u = (a/b)^{1/({m-1})}$ yields $v^{\prime}(t) + a v(t) = amv(\alpha t),$ for which Lemma~1 yields the critical values $\alpha_n = m^{1/n}$. We expect the remaining steps, including the moment matching, to proceed as in the $m = 2$ case.} \red{Yet} another avenue of investigation is how to solve~(\ref{eq:nonlinear}) for other initial conditions.   One possibility might be to 
use the $u(0) = 1$ solutions in combination with continuation on the $u_0$ parameter 
to generate solutions corresponding to nearby initial conditions, \red{and hence fully explore the $(u(0), \alpha)$ parameter space}.

\section{Acknowledgments}
{The authors thank the anonymous referee\red{s} for valuable feedback.} The second author thanks the staff, faculty, and colleagues of the Department of Mathematical Sciences at Stellenbosch University for creating a welcoming and productive environment that made this research possible while on sabbatical leave.

\section{Declarations}
\textbf{Ethical approval:} Not applicable.\\
\textbf{Competing interests:} The authors declared no potential conflicts of interest concerning this article’s research, authorship, and publication.\\
\textbf{Generative AI:} \red{Generative AI (ChatGPT 5.5) was used to assist in editing the authors’ text for spelling, grammar, and structure in this work. AI was not used in the generation of any code, numerical results, or figures.}


\appendix%
\section{Appendix}\label{sec:appendix}%
The details of computing the scaling factors, $C_n$, given in Table~\ref{table:c}
are as follows.   Denote the $j$-th moments of $v$ by
\be
\mu_j(v) = \int_0^\infty t^j v(t) \, dt, \quad j = 0, 1, 2, \ldots
\ee
{Under our assumption that $v(t) = {{O}(e^{-t})}$ as $t\to\infty$, all the $\mu_j(v)$ are bounded.}
By multiplying~(\ref{eq:v}) by $t^j$ and integrating by parts, it follows that\footnote{Here, and below, moments with negative indices are included for ease of notation. They should be taken to be identically zero.}
\be
(\alpha^{j+1} - 2) \mu_j(v) - j \alpha^{j+1} \mu_{j-1}(v)
- \mu_j(v^2) = 0,
\label{eq:moments}
\ee
where the facts that $v(0)=0$ and 
$t^j v(t) \to 0$   as $t \to \infty$ were used  to eliminate the derivative term. 

Now, for $n = 1, 2, \ldots$, consider the perturbation $\alpha = \alpha_{n}+\epsilon$, where $\alpha_n = 2^{1/n}$, and the expansion $v(t) = \epsilon v_0(t) + \epsilon^2 v_1(t) + 
\cdots$, as in section~\ref{sec:nonlinear}. As before,  $v_0$ is the solution to~(\ref{eq:linear2}) with $\alpha= \alpha_n$.  
Substitute this into~(\ref{eq:moments}) and use the
binomial theorem to expand $(\alpha_{n}+\epsilon)^{j+1}$.
Equating the first two powers of $\epsilon$ gives
\be  
\label{eqn:jOd}
(\alpha^{j+1}_{n}-2)\mu_j(v_0) -j \alpha^{j+1}_{n}\mu_{j-1}(v_0) = 0,
\ee
and  
\be  \label{eqn:jOd2a}
(\alpha^{j+1}_{n}-2)\mu_j(v_1)-j \alpha^{j+1}_{n}\mu_{j-1}(v_1)
+(j+1)\alpha_n^j\mu_{j}(v_0) - j(j+1)\alpha_n^j\mu_{j-1}(v_0) - \mu_j(v_0^2) = 0,
\ee
respectively. Using the former expression, the $\mu_{j-1}(v_0)$ term in the latter can be eliminated, giving
\be  \label{eqn:jOd2}
(\alpha^{j+1}_{n}-2)\mu_j(v_1)-j \alpha^{j+1}_{n}\mu_{j-1}(v_1)
+ \frac{2(j+1)}{\alpha_{n}}\mu_j(v_0) - \mu_j(v_0^2) = 0.
\ee

At this stage it is useful to note that,  because $v_0$ is known---at least, up to the scaling constants $C_n$ that we are presently trying to determine---all moments involving 
 $v_0$ can be computed numerically; see~(\ref{eq:summations}) below.   This is not the case for moments of $v_1$, however, since $v_1$ is unknown, and hence the purpose in the following is to eliminate such moments.

With $j = n-1$, we have $\alpha_{n}^{j+1} = 2$ and equations~(\ref{eqn:jOd}) and~(\ref{eqn:jOd2}) simplify to 
\be\label{eqn:nOd}  \mu_{n-2}(v_0) = 0, \ee
and 
\be\label{eqn:nOd2}  
 -2(n-1)\mu_{n-2}(v_1)
+ \frac{2n}{\alpha_{n}}\mu_{n-1}(v_0) - \mu_{n-1}(v_0^2) = 0,
\ee
respectively.

The fact that the $\mu_{n-1}(v_1)$ term is eliminated in the latter expression is precisely the reason for considering the $(n-1)$th moment. This  will allow us to develop a telescoping recurrence in order to eliminate moments involving the unknown $v_1$. To this end, note that~(\ref{eqn:nOd2})
can be rearranged as
\be\label{eqn:nOd2_rearranged} 
\frac{2n}{\alpha_{n}}\mu_{n-1}(v_0) = \mu_{n-1}(v_0^2) + 2(n-1)\mu_{n-2}(v_1),\ee
 which we return to momentarily. 

Now consider $j < n-1$. It follows   from~(\ref{eqn:jOd}) 
and~(\ref{eqn:nOd})  that 
\be  \mu_{j}(v_0) = 0, \qquad j < n-1.\ee
This causes the third term
on the left of~(\ref{eqn:jOd2}) to vanish, which gives the
recurrence
\be\label{eq:mu} \mu_j(v_1) = \frac{1}{\alpha_{{n}}^{j+1}-2}\left(\mu_j(v_0^2) + j\alpha_{{n}}^{j+1}\mu_{j-1}(v_1)\right), \qquad j < n-1.\ee

Returning now to~(\ref{eqn:nOd2_rearranged}) and   applying~(\ref{eq:mu}) with $j=n-2$ gives
\begin{eqnarray}
\frac{2n}{\alpha_{n}}\mu_{n-1}(v_0) = \mu_{n-1}(v_0^2) + \frac{2(n-1)}{\alpha^{n-1}_{n}-2}\left(\mu_{n-2}(v_0^2) + (n-2)\alpha_{n}^{n-1}\mu_{n-3}(v_1)\right). 
\label{eq:recur}
 \end{eqnarray}

If $n=2$ the term on the right involving $v_1$ is eliminated and, with 
 $\alpha_2 =   2^{1/2}$, one then gets
 \be
2\sqrt{2} \,  \mu_1(v_0) = 
\mu_1(v_0^2) - (2+\sqrt{2}) \, \mu_0(v_0^2).
\label{eq:neq1}
 \ee
 If $n > 2$ then further substitutions of~(\ref{eq:mu}) into (\ref{eq:recur}) are made until the resulting series likewise terminates with the vanishing of the quantity involving $v_1$ and only moments of $v_0$ and $v_0^2$ remain.  
 The general formula involves the following linear combination of the $\mu_j(v_0^2)$ on the right
 \be\label{eq:gen1}
\frac{2n}{\alpha_{n}} \mu_{n-1}(v_0)
= \sum_{j=0}^{n-1} w_j \mu_j(v_0^2), \qquad n = 1, 2, \ldots,
 \ee
 where the $w_j$ are generated in the reverse direction by
 $w_{n-1} = 1$, and
\be
 w_{j-1} = \frac{j\alpha^{j+1}_{n}}{\alpha^{j}_{n}-2}w_{j}, \quad  j = n-1, \ldots, 1.
 \label{eq:gen}
 \ee

Letting $v_0(t) = C_nE(t;\alpha_{n})$ as in~(\ref{eq:series}), it follows that~(\ref{eq:gen1}) becomes 
 \be\label{eqn:gens} \frac{2n}{\alpha_{n}}C_n\mu_{n-1}(E_{n}) = C_n^2 \sum_{j=0}^{n-1}w_j\mu_{j}(E^2_{n}), \ee
 where we have denoted $E(t; \alpha_{n})$ by $E_{n}$ for brevity. Finally, 
solving for $C_n$  gives 
\be\label{eq:scalingcoeffs}
C_{n} = \frac{2n\mu_{n-1}(E_{n}) }{\alpha_{n} \displaystyle\sum_{j=0}^{n-1}w_j\mu_{j}(E^2_{n})}, \quad n = 1, 2, \ldots.\ee

For $n = 1$ ($\alpha_1 = 2$) one
gets
\be\label{eq:scalingcoeffs1} C_1  = \frac{\mu_0(E_1)}{\mu_0(E_1^2)},\ee
a formula already 
given in~(\ref{eq:C1}).
For $n = 2$ ($\alpha_2 = \sqrt{2}$)   
the formula is 
\be\label{eq:scalingcoeffs2} C_2 = 
\frac{2\sqrt{2} \,\mu_1(E_{2})}
{\mu_1(E^2_2)-(2+\sqrt{2})\mu_0(E_2^2)},\ee
 which also follows directly from~(\ref{eq:neq1}).
For larger values of $n$ the formulas are not as clean as the surds do not cancel elegantly.

Explicit series formulas for the moments in~(\ref{eq:scalingcoeffs})--(\ref{eq:scalingcoeffs2}) can be obtained by integrating the series~(\ref{eq:series}) term-by-term. In particular, 
\be
\mu_j(E_n) = j!\sum_{k=0}^\infty \frac{b_k}{\alpha_n^{(j+1)k}}
\qquad \text{and} \qquad 
\mu_j(E^2_n) = j!\sum_{k=0}^\infty\sum_{\ell=0}^\infty \frac{b_k b_\ell}{(\alpha_n^k+\alpha_n^\ell)^{j+1}},
\label{eq:summations}
\ee
where 
\be
b_k = \left\{\begin{array}{cl}
1, & k = 0,\\
\displaystyle\frac{2^k}{(1-\alpha_n)\cdots(1-\alpha_n^k)}, & 
k>0.  
\end{array}\right.
\ee
Although the coefficients in these series have exponential asymptotic decay, they have alternating signs and large transient growth, making them numerically unstable, particularly when $n$ is large. Extended precision is therefore required to compute the moments accurately from these formulas. {Alternatively, one could approximate the moments via numerical integration after solving~(2.1) using the spectral method.}

\bibliographystyle{siamplain} 
\bibliography{mybib}%

\end{document}